\documentclass[twocolumn]{IEEEtran}
\usepackage{amsmath,amssymb,amsthm,epsfig,color,subfigure,empheq,graphicx,graphics,balance}
\usepackage{enumerate,url,algorithm,algorithmic,wasysym,epstopdf,enumitem}
\usepackage{multirow, nicefrac} 
\usepackage{fancyvrb}
\VerbatimFootnotes
\usepackage{xcolor}
\usepackage{amsfonts}
\usepackage{xparse}
\newsavebox{\fminipagebox}
\NewDocumentEnvironment{fminipage}{m O{\fboxsep}}
 {\par\kern#2\noindent 
  \vspace{-3ex}
 \begin{center}\begin{lrbox}{\fminipagebox}
  \begin{minipage}{#1}\ignorespaces}
 {\end{minipage}\end{lrbox}%
  \makebox[#1]{%
    \kern\dimexpr-\fboxsep-\fboxrule\relax
    \fbox{\usebox{\fminipagebox}}%
    \kern\dimexpr-\fboxsep-\fboxrule\relax
  }\par\kern#2 \end{center}
  \vspace{-1ex}
 }

\newcounter{ProblemCounter}
\renewcommand\theProblemCounter{\arabic{ProblemCounter}}

\DeclareMathOperator{\range}{range}

\DeclareMathOperator{\nullspace}{null}

\DeclareMathOperator{\trace}{Tr}

\newtheorem{assumption}{Assumption}

\newtheorem{lemma}{Lemma}

\newtheorem{theorem}{Theorem}

\newtheorem{remark}{Remark}

\newcommand \bzero{\mathbf{0}}

\newcommand \ba{\mathbf{a}}
\newcommand \bb{\mathbf{b}}

\newcommand \bd{\mathbf{d}}
\newcommand \be{\mathbf{e}}
\newcommand \bg{\mathbf{g}}
\newcommand \bh{\mathbf{h}}
\newcommand \bi{\mathbf{i}}

\newcommand \bn{\mathbf{n}}

\newcommand \bp{\mathbf{p}}
\newcommand \bq{\mathbf{q}}

\newcommand \bu{\mathbf{u}}
\newcommand \bv{\mathbf{v}}
\newcommand \bw{\mathbf{w}}
\newcommand \bx{\mathbf{x}}

\newcommand \bz{\mathbf{z}}
\newcommand \bA{\mathbf{A}}
\newcommand \bB{\mathbf{B}}

\newcommand \bD{\mathbf{D}}

\newcommand \bG{\mathbf{G}}

\newcommand \bI{\mathbf{I}}
\newcommand \bJ{\mathbf{J}}
\newcommand \hbJ{\hat{\mathbf{J}}}

\newcommand \bL{\mathbf{L}}
\newcommand \bM{\mathbf{M}}

\newcommand \bS{\mathbf{S}}

\newcommand \bU{\mathbf{U}}
\newcommand \bV{\mathbf{V}}

\newcommand \bY{\mathbf{Y}}
\newcommand \bZ{\mathbf{Z}}

\newcommand \bgamma{\boldsymbol{\gamma}}

\newcommand \bepsilon{\boldsymbol{\epsilon}}

\newcommand \btheta{\boldsymbol{\theta}}

\newcommand \blambda{\boldsymbol{\lambda}}
\newcommand \bmu{\boldsymbol{\mu}}

\newcommand \mcA{\mathcal{A}}

\newcommand \mcD{\mathcal{D}}
\newcommand \mcE{\mathcal{E}}

\newcommand \mcG{\mathcal{G}}

\newcommand \mcL{\mathcal{L}}

\newcommand \mcN{\mathcal{N}}

\newcommand \mcZ{\mathcal{Z}}

\newcommand \cbx{\check{\mathbf{x}}}

\newcommand \hbx{\hat{\mathbf{x}}}

\newcommand \bbgamma{\bar{\boldsymbol{\gamma}}}

\renewcommand{\d}[1]{\ensuremath{\operatorname{d}\!{#1}}}

\begin{document}
	\begin{table*}
		\Large
\copyright 2021 IEEE.  Personal use of this material is permitted.  Permission from IEEE must be obtained for all other uses, in any current or future media, including reprinting/republishing this material for advertising or promotional purposes, creating new collective works, for resale or redistribution to servers or lists, or reuse of any copyrighted component of this work in other works.
	\end{table*}
	
\pagebreak
\title{Learning to Solve the AC-OPF using Sensitivity-Informed Deep Neural Networks}

\author{Manish K. Singh~\IEEEmembership{Member,~IEEE},~Vassilis Kekatos~\IEEEmembership{Senior Member,~IEEE},\\ and Georgios~B.~Giannakis~\IEEEmembership{Fellow,~IEEE}}


\maketitle
\begin{abstract}
To shift the computational burden from real-time to offline in delay-critical power systems applications, recent works entertain the idea of using a deep neural network (DNN) to predict the solutions of the AC optimal power flow (AC-OPF) once presented load demands. As network topologies may change, training this DNN in a sample-efficient manner becomes a necessity. To improve data efficiency, this work utilizes the fact OPF data are not simple training labels, but constitute the solutions of a parametric optimization problem. We thus advocate training a sensitivity-informed DNN (SI-DNN) to match not only the OPF optimizers, but also their partial derivatives with respect to the OPF parameters (loads). It is shown that the required Jacobian matrices do exist under mild conditions, and can be readily computed from the related primal/dual solutions. The proposed SI-DNN is compatible with a broad range of OPF solvers, including a non-convex quadratically constrained quadratic program (QCQP), its semidefinite program (SDP) relaxation, and MATPOWER; while SI-DNN can be seamlessly integrated in other learning-to-OPF schemes. Numerical tests on three benchmark power systems corroborate the advanced generalization and constraint satisfaction capabilities for the OPF solutions predicted by an SI-DNN over a conventionally trained DNN, especially in low-data setups.
\end{abstract}


\section{Introduction}\label{sec:intro}
\allowdisplaybreaks
Power system operation involves routinely solving various renditions of the optimal power flow (OPF) task. Physical expansion of power networks, increasing number of dispatchable resources, and renewable generation-induced volatility call for solving large-scale OPF problems frequently. These problems naturally involve nonlinear alternating-current (AC) power flow equations as constraints. Such formulations, referred to as AC-OPF, are non-convex and are often computationally too expensive for real-time applications. Approximate formulations such as the linearized or so termed DC-OPF serve as the pragmatic resort for such setups.

A concerted effort towards handling AC-OPF efficiently has yielded popular nonlinear solvers such as MATPOWER~\cite{MATPOWER}, and efficient conic relaxations with optimality guarantees for frequently encountered problem instances~\cite{Bai08}, \cite{Bose}, \cite{MSL15}, \cite{Low14}. Despite significant advancements in numerical solvers, scalability of the AC-OPF may still be a challenge, particularly in online, combinatorial, and stochastic settings~\cite{Xavier20SCUC}, \cite{Fioretto1}. To alleviate these issues, there has been growing interest in developing machine learning-based approaches (particularly deep learning) for OPF~\cite{Xavier20SCUC}, \cite{Fioretto1}, \cite{DeepOPFPan19}, \cite{Pan20DeepOPFplus}, \cite{pan2020feasopt}, \cite{GuhaACOPF}, \cite{ZamzamBaker19}; see~\cite{YueZhao_DL4P} for recent applications. The primary advantage of machine learning-based approaches lies in the speed-up during the inference phase. For instance, compared to conventional solvers, deep learning-based approaches have offered speed-ups by factors as high as 200 for DC-OPF, and 35 for AC-OPF~\cite{Pan20DeepOPFplus}, \cite{pan2020feasopt}. 

To entirely bypass numerical solvers for the OPF, one can opt for unsupervised learning approaches, such as the ones suggested in~\cite{RibeiroGNNOPF}, \cite{GKJ2020}, \cite{lange2020learning}, \cite{GMDK21}, \cite{OPFandLearnTSG21}. However, the performance of OPF solvers for offline data generation and the availability of historical data with power utilities adequately motivate supervised learning. Nevertheless, there are two main challenges central to learning for OPF. First, traditional deep learning approaches are not amenable to enforcing constraints even for the training set. Predictions for OPF minimizers may have limited standalone value if the related constraints are violated. Second, power systems undergo frequent topological and operational changes as generating units can be (de)-committed and transmission lines or bus/line reactors can be switched. Such changes may require retraining a DNN potentially at short intervals~\cite{chen2020metalearning}. Deep learning-based approaches are traditionally data-intense and frequent retraining for large systems may be prohibitive.

To cope with the first challenge, a DNN may be engaged to better initialize existing numerical solvers~\cite{ZamzamBaker19}, or to predict active constraints and thus result in an OPF model with much reduced number of constraints~\cite{GuhaACOPF},~\cite{ChenZhang20},~\cite{DekaMisraPowerTech19}. Another group of approaches targets constraint satisfaction by penalizing constraint violations and the related Karush–Kuhn–Tucker (KKT) conditions~\cite{DeepOPFPan19}, \cite{pan2020feasopt}, or explicitly resorting to a Lagrangian dual scheme for DNN training~\cite{Fioretto1}, \cite{GKJ2020}, \cite{Nandwani19}. The third alternative involves post-processing DNN predictions by projecting them using a power flow solver~~\cite{DeepOPFPan19},~\cite{ZamzamBaker19}. Although the projected point satisfies the power flow equations, it may still violate engineering limits. 

To cater to the second challenge of frequent model changes, sample-efficient learning models that generalize well are well motivated. One way to achieve that is via meta-learning, according to which training datasets generated from diverse grid topologies are used to train a vector of meta-weights~\cite{chen2020metalearning}. When it later comes to training a DNN to handle a particular topology, its weights are initialized to the values of meta-weights. Alternatively, the sample efficiency of learning models could be improved by prudently designing the DNN architecture upon leveraging prior information. For example, by seeking an input-convex DNN when the underlying input-output mapping is convex~\cite{Zhang20icnn}, or using DNNs that \emph{unroll} an iterative optimization algorithm~\cite{ZhangWangGiannakis19}, or adopting graph-based priors~\cite{RibeiroGNNOPF}, \cite{yang2020robust}. The previously mentioned approaches that design objective functions based on OPF constraints and KKT conditions may also be seen as including prior knowledge in training, thus enhancing learnability. 

Recent works in the general area of \emph{physics-informed learning} aim at incorporating prior information on the underlying data, not occurring in conventional training datasets~\cite{Raissi19}. Specifically, for learning solutions of differential equations, the derivatives of DNN output with respect to input carry obvious virtue. Suitably utilizing these derivatives yields significant advantages in such applications~\cite{Raissi19},~\cite{misyris2020PINN}. While for differential equations derivatives with respect to time and/or space dimensions emerge naturally, sensitivities in an optimization problem have been underappreciated. To this end, our recent work proposed a novel approach of training sensitivity informed DNNs (SI-DNN), intended to learn OPF solutions~\cite{SGKCB2020}. Training SI-DNNs requires computing the sensitivities of OPF minimizers with respect to the input parameters. For DC-OPF posed as linear or quadratic programs, computing these sensitivities is simpler, and using these to train DNNs yielded remarkable improvements; see~\cite{SGKCB2020}. 

\emph{Contributions:} The contributions of this work are on four fronts: \emph{c1)} We put forth a novel approach for training DNNs to predict AC-OPF solutions by matching not only the OPF minimizers, but also their sensitivities (partial derivatives collected in a Jacobian matrix) with respect to the OPF problem parameters (e.g., load demands). \emph{c2)} We compute the desired sensitivities for general nonlinear OPF formulations building upon classical works on perturbation analysis of continuous optimization problems~\cite{shapiro2013perturbation}, \cite{fiacco1976sensitivity}, \cite{Conejo06}. Existing sensitivity results assume certain constraint qualifications that may not be satisfied by OPF instances. \emph{c3)} We relax such assumptions and establish that the sensitivities of primal OPF solutions do exist under milder conditions. \emph{c4)} In pursuit of globally optimal OPF solutions for training a DNN, we also study the SDP formulation of the AC-OPF and compute its sensitivities by utilizing the sensitivities of the related QCQP. Such shortcut obviates the difficulty of differentiating a conic program. It is worth stressing that the proposed sensitivity-informed methodology can be used in tandem with other learning approaches, such as those reviewed earlier imposing physics-inspired DNN architectures, penalizing constraint violations, or adopting meta-learning.

\emph{Motivation:} The advantage of using a DNN to predict OPF minimizers is its computational speed. Running a DNN is significantly faster than running an OPF solver. Alternatively, the DNN prediction can be used to \emph{warm-start} an OPF solver and thus expedite its convergence. This speed-up advantage during testing (inference) is featured by any DNN. Sensitivity-informed DNNs are advocated here as a means to expedite training as well. This is possible because an SI-DNN can attain the same prediction accuracy as a plain DNN using much fewer samples. Our numerical tests demonstrate that depending on the network an SI-DNN requires 1/4 or even 1/10 of the training samples compared to a plain DNN. Such speed-up is significant in setups where there is not sufficient time for generating samples offline. Such setups could arise when a DNN is used to predict the setpoints of inverter-interfaced distributed energy resources and the feeder is reconfigured frequently; or when unit commitment decisions change regularly the set of online generators in a transmission system; or when the DNN is part of a stochastic unit commitment formulation and has to be retrained for the various commitment configurations visited by a branch-and-bound algorithm.

The rest of the paper is organized as follows: Section~\ref{sec:SIT} presents the key methodology for SI-DNNs. Section~\ref{sec:QCQP} poses the AC-OPF as a parametric optimization taking the form of a non-convex QCQP. Section~\ref{sec:sensitivity} establishes that the sought sensivities of the AC-OPF exist under mild conditions and explains how they can be readily computed. Section~\ref{sec:SDP} computes the sensitivities of globally optimal AC-OPF solutions obtained via the SDP relaxation of the OPF. The numerical tests of Section~\ref{sec:tests} contrast SI-DNN to conventionally trained DNNs using datasets generated by MATPOWER and the SDP-based OPF solver on three benchmark power systems. Conclusions are drawn in Section~\ref{sec:conclusions}.

\emph{Notation}: lower- (upper-) case boldface letters denote column vectors (matrices). Calligraphic symbols are reserved for sets. Symbol $^{\top}$ stands for transposition, vectors $\mathbf{0}$ and $\mathbf{1}$ are the all-zeros and all-ones vectors or matrices, and $\be_n$ is the $n$-th canonical vector of appropriate dimensions implied by the context. Operator $\mcD(\bx)$ returns a diagonal matrix with the entries of vector $\bx$ on its main diagonal.

\section{Sensitivity-Informed Training for DNNs}\label{sec:SIT}
Consider the task of optimally dispatching generators and flexible loads to meet power balance constraints $\bh(\cdot)$ while enforcing engineering limits captured by $\bg(\cdot)$. OPF boils down to a \emph{parametric optimization problem} that has to be solved routinely for different values of the problem parameters being the time-varying renewable generation and loads. The problem can be abstracted as: Given a vector $\btheta\in\mathbb{R}^P$ of problem parameters or inputs, find an optimal dispatch $\bx_{\btheta}\in\mathbb{R}^N$ consisting of voltages and generator setpoints as the minimizer
\begin{align}\label{eq:Ptheta}
\bx_{\btheta}\in\arg\min_{\bx}~&f(\bx;\btheta)\tag{$P_{\btheta}$}\\
\mathrm{s.to}~&~\bh(\bx;\btheta)= \bzero:~\blambda_{\btheta}\nonumber\\
~&~\bg(\bx;\btheta)\leq \bzero:~\bmu_{\btheta}\nonumber
\end{align}
where functions $f(\bx;\btheta)$, $\bh(\bx;\btheta)$, and $\bg(\bx;\btheta)$ are continuously differentiable with respect to $\bx$ and $\btheta$; and vectors $(\blambda_{\btheta},\bmu_{\btheta})$ collect the optimal dual variables corresponding to the equality and inequality constraints, respectively.

To save on running time and computational resources, rather than solving \eqref{eq:Ptheta}, one can adopt a \emph{learning-to-optimize approach} and train a learning model such as a DNN to predict approximate solutions of \eqref{eq:Ptheta}; see e.g.,~\cite{Sidiropoulos18}. Once presented with an instance of $\btheta$, this DNN can be trained to output a predictor $\hbx(\btheta;\bw)$ of $\bx_{\btheta}$. The DNN is parameterized by weights $\bw$, which can be selected upon minimizing a suitable distance metric or loss function between $\bx_{\btheta}$ and $\hbx(\btheta;\bw)$ over a training set. 

Given a choice for a DNN architecture, the straightforward approach for learning-to-optimize entails two steps:
\renewcommand{\labelenumi}{\emph{S\arabic{enumi})}}
\begin{enumerate}
    \item Building a labeled training dataset $\{\btheta_s,\bx_{\btheta_s}\}_{s=1}^S$ by solving $S$ instances of \eqref{eq:Ptheta}; and
    \item Learning $\bw$ by minimizing a data fitting loss over the training dataset as
\begin{equation}\label{eq:p-dnn}
\min_{\bw}\sum_{s=1}^S \ell\left(\hbx(\btheta_s;\bw),\bx_{\btheta_s}\right).
\end{equation}
\end{enumerate}

For a regression task such as the one considered here, commonly used loss functions $\ell$ include the mean squared error (MSE) $\|\hbx(\btheta;\bw)-\bx_{\btheta}\|_2^2$, or the mean absolute error (MAE) $\|\hbx(\btheta;\bw)-\bx_{\btheta}\|_1$. We refer to a DNN trained by solving \eqref{eq:p-dnn} as a \emph{plain DNN or P-DNN} for short. The conventional P-DNN approach focuses merely on the dataset $\{\btheta_s, \bx_{\btheta_s}\}_{s=1}^S$, and is oblivious of any additional properties the mapping $\btheta\rightarrow\bx_{\btheta}$ induced by \eqref{eq:Ptheta} bears. 

The key idea here is to extend each training data pair $(\btheta_s,\bx_{\btheta_s})$ as $(\btheta_s,\bx_{\btheta_s},\bJ_{\btheta_s})$, where $\bJ_{\btheta_s}:=\left[\nabla_{\btheta}\bx_{\btheta}\right]_{\btheta=\btheta_s}$ is the Jacobian matrix carrying the partial derivatives of the minimizer $\bx_{\btheta_s}$ with respect to $\btheta$ evaluated at $\btheta=\btheta_s$ assuming such sensitivities exist. To incorporate the sensitivity information into DNN training, we propose augmenting the loss function with an additional fitting term as
\begin{equation}\label{eq:si-dnn}
\min_{\bw}\sum_{s=1}^S\|\hat{\bx}(\btheta_s;\bw)-\bx_{\btheta_s}\|_2^2 +\rho \|\hbJ(\btheta_s;\bw)-\bJ_{\btheta_s}\|_F^2
\end{equation}
where $\rho>0$ is a scalar weight and $\|\cdot\|_F$ denotes the matrix Frobenius norm. The DNN trained by solving \eqref{eq:si-dnn} aims to match not only the target output $\bx_{\btheta}$, but also the sensitivities of $\bx_{\btheta}$ with respect to $\btheta$. We term this neural network a \emph{sensitivity-informed DNN or SI-DNN}. 

\begin{figure*}[t]
    \centering
    \hspace*{-1em}
    \includegraphics[scale=0.66]{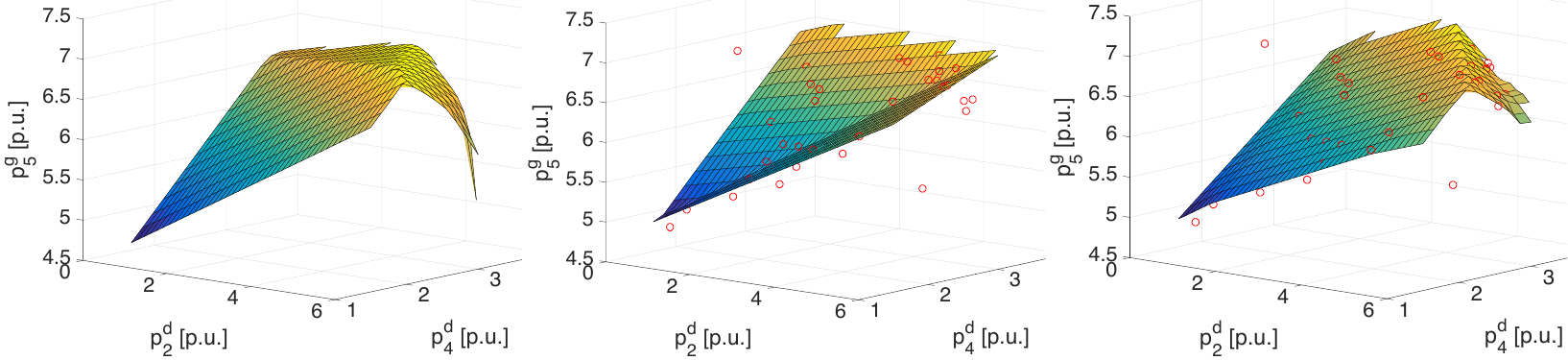} 
    \caption{The left panel depicts the optimal generation dispatch $p^g_5(\btheta)$ for bus 5 as a function of load demands at buses 2 and 4, that is $\btheta:=[p^d_2~p^d_4]^\top$. Sampling the parameter space of $\btheta$'s provided 523 feasible OPF instances, of which 37 instances constituted the training set. The center and right panel show the dispatches learned by a P-DNN and an SI-DNN. For the two DNNs, we used the same training points (red circles), architecture (two hidden layers with 16 neurons each), optimizer, and learning rates. P-DNN fails to learn the drop in $p^g_5$ for larger load demands.}
    \label{fig:toy}
\end{figure*}

%
%

To better understand \eqref{eq:si-dnn} and motivate the inclusion of sensitivities, let us put forth the ensuing interpretation. Consider learning function $x:\mathbb{R}\rightarrow\mathbb{R}$, which can be an OPF mapping $x(\theta)$ with $P=N=1$. Under the typical learning setup, one aims to build a DNN $\hat{x}(\theta)$ and approximate $x(\theta)$ given training samples $\{(\theta_s,x(\theta_s))\}_{s=1}^S$. Different from this setup, suppose we are given the additional information of function derivative values at the training samples, so that the training dataset consists of the triplets $\{(\theta_s,x(\theta_s),x'(\theta_s))\}_{s=1}^S$. The pertinent question is how to utilize the extra sensitivity information.

Linearizing function $x$ around a sample $\theta_s$ yields \[x(\theta_s+\epsilon)\simeq x(\theta_s)+\epsilon \cdot x'(\theta_s)\] 
for any small $\epsilon$. Linearizing the DNN output yields similarly
\[\hat{x}(\theta_s+\epsilon)\simeq \hat{x}(\theta_s)+\epsilon\cdot \hat{x}'(\theta_s)\]
where $\hat{x}'(\theta_s)$ is the derivative of the DNN output with respect to its input $\theta$ evaluated at $\theta_s$. Suppose now $\epsilon$ is a zero-mean random variable with variance $\mathbb{E}[\epsilon^2]=\sigma^2$. Instead of training the DNN by minimizing the loss $(\hat{x}(\theta_s)-x(\theta_s))^2$ summed up over all $s$, one can aim at fitting the function around a sphere of essential radius $\sigma$ that is centered at $\theta_s$ by minimizing
\begin{align*}
&\mathbb{E}_\epsilon\left[\left(\hat{x}(\theta_s+\epsilon)-x(\theta_s+\epsilon)\right)^2\right] \\
&\simeq \left(\hat{x}(\theta_s)-x(\theta_s)\right)^2+\sigma^2\left(\hat{x}'(\theta_s)-x'(\theta_s)\right)^2.
\end{align*}
Of course, the previous loss is also summed up over all $s$. Interestingly, this stochastic interpretation of function fitting yields the sensitivity-aware training task of \eqref{eq:si-dnn} upon identifying $\rho=\sigma^2$. This scalar case of $P=N=1$ can be trivially extended to a general (OPF) mapping of arbitrary dimensions $P$ and $N$ upon substituting $\epsilon$ by a zero-mean random vector $\bepsilon\in\mathbb{R}^P$ with covariance matrix $\mathbb{E}[\bepsilon\bepsilon^\top]=\sigma^2\bI_P$ and replacing derivatives with Jacobian matrices. This interpretation not only justifies the form of \eqref{eq:si-dnn}, but also explains geometrically how sensitivity-informed training uses the point information $\{(\theta_s,x(\theta_s),x'(\theta_s))\}_{s=1}^S$ to extrapolate in a neighborhood around each training datum.  


SI-DNNs for learning optimizers were first introduced for solving multiparametric QPs (MPQP) in the conference precursor of this work~\cite{SGKCB2020}. For MPQPs, the minimizer $\bx_{\btheta}$ is known to be a piecewise affine function of $\btheta$~\cite{BBM03}, \cite{TJKT20}. Hence, a DNN with rectified linear unit (ReLU) activations is well-motivated as it can describe such mapping. If hypothetically trained to zero training error, this SI-DNN would yield perfectly accurate predictions in a neighborhood of each training datum $\btheta_s$. Numerical tests showed improvements of 2-3 orders of magnitude for SI-DNN over P-DNN in inferring MPQP solutions~\cite{SGKCB2020}.

This work advocates that the sample efficiency benefit of SI-DNN over P-DNN goes well beyond MPQPs. Before delving into the details and for the sake of visualization, we present some numerical tests on a toy 5-bus power system; a proper numerical evaluation of SI-DNN is deferred to Section~\ref{sec:tests}. This PJM 5-bus system was dispatched via AC-OPF by varying the active load demands $\btheta:=[p^d_2~p^d_4]^\top$ on buses 2 and 4 within [1.5, 3.75] per unit (pu) and [0.4, 6] pu, respectively. Fixing the other demands, we dispatched the generators at buses 1 and 5, and removed other generators. Figure~\ref{fig:toy} depicts the performance improvement of SI-DNN over P-DNN in predicting $p^g_5$, the optimal dispatch at bus 5.

It is worth clarifying that this work does not train a DNN to predict OPF solutions under \emph{different} power system topologies. On the contrary, it aims at learning the OPF mapping for a single given topology under diverse loading conditions. The fact that the network topology may be changing across time justifies the need to improve on data efficiency, so that after a topology change, the corresponding DNN can be trained afresh using fewer OPF examples generated using the new topology.

This new learning-to-optimize approach alters the two steps of the P-DNN workflow as follows: For step \emph{S1)}, in addition to the minimizer $\bx_{\btheta_s}$, we now have to compute the Jacobian matrix $\bJ_{\btheta_s}$ for all instances $s$, if such sensitivities exist. Sections~\ref{sec:sensitivity} and \ref{sec:SDP} explain how and when such Jacobian matrices can be computed for a non-convex and a convexified rendition of the AC-OPF. The punchline is that obtaining $\bJ_{\btheta}$ requires minimal additional computational effort and no intervention to the OPF solver. Once the primal/dual solutions have been found by the OPF, computing $\bJ_{\btheta}$ is as simple as solving a linear system of equations. 

For step \emph{S2)}, we migrate from solving \eqref{eq:p-dnn} to \eqref{eq:si-dnn}. Matrix $\bJ_{\btheta_s}$ is a constant that has been evaluated for each $s$ during \emph{S1)}. Matrix $\hbJ(\btheta_s;\bw)$ on the other hand is a function of $\bw$ and is not straightforward to compute. Fortunately, computing $\hbJ(\btheta_s;\bw)$ can be performed efficiently thanks to advances in automatic differentiation~\cite{autodiff}. Modeling \eqref{eq:si-dnn} in existing DNN software platforms (e.g., TensorFlow) is almost as easy as modeling \eqref{eq:p-dnn} modulo the coding modifications deferred to Appendix~\ref{sec:AppA}. Our tests of Section~\ref{sec:tests} further show that the extra computational time for solving \eqref{eq:si-dnn} is modest. 
Before computing $\bJ_{\btheta}$, we first pose AC-OPF as a parametric QCQP.

\section{AC-OPF as a Parametric QCQP}\label{sec:QCQP}
A power network with $N_b$ buses can be represented by an undirected connected graph $\mcG=(\mcN,\mcE)$, whose nodes $n\in\mcN:= \{1,\ldots,N_b\}$ correspond to buses, and edges $e=(n,k)\in\mcE$ to transmission lines, with cardinality $|\mcE|:=E$. Given line impedances, one can derive the $N_b\times N_b$ bus admittance matrix $\bY=\bG+j\bB$. Let $v_n=v^r_n+jv^i_n$ and $p_n+jq_n$ denote respectively the complex voltage and power injection at bus $n\in\mcN$. Power injections are quadratically related to bus voltages through the power flow equations
\begin{align*}
p_n&=\sum_{k=1}^{N_b} v^r_n(v^r_kG_{nk} - v^i_kB_{nk}) +  v^i_n(v^i_kG_{nk} + v^r_kB_{nk})\\
q_n&=\sum_{k =1}^{N_b} v^i_n(v^r_kG_{nk} - v^i_kB_{nk}) - v^r_n(v^i_kG_{nk} + v^r_kB_{nk}).
\end{align*}
If $\bv\in\mathbb{R}^{2N_b}$ collects the real and imaginary parts of nodal voltages as $\bv:=[(\bv^r)^\top~(\bv^i)^\top]^\top$ with $\bv^r:=\{v^r_n\}_{n=1}^{N_b}$ and $\bv^i:=\{v^i_n\}_{n=1}^{N_b}$, the power flow equations can be written as
\begin{subequations}\label{eq:pfc}
\begin{align}
p_n&=\bv^{\top}\bM_{p_n}\bv\label{eq:pfc:p}\\
q_n&=\bv^{\top}\bM_{q_n}\bv\label{eq:pfc:q}
\end{align}
\end{subequations}
where $\bM_{p_n}$ and $\bM_{q_n}$ are $2N_b\times2N_b$ symmetric real-valued matrices~\cite{redux}. Squared voltage magnitudes can also be expressed as quadratic functions of $\bv$ as
\begin{align}\label{eq:volts}
|v_n^r+jv_n^i|^2=\bv^{\top}\bM_{v_n}\bv
\end{align}
where $\bM_{v_n}:=\be_n\be_n^\top+\be_{N_b+n}\be_{N_b+n}^\top$. The same holds true for the squared magnitude of line currents. If $y_{mn}$ is the series admittance of line $(m,n)\in\mcE$, the current flowing on this line is $\tilde{i}_{mn}=(\tilde{v}_m-\tilde{v}_n)y_{mn}$, and thus,
\begin{equation}
|\tilde{i}_{mn}|^2=\bv^{\top}\bM_{i_{mn}}\bv
\end{equation}
where $\bM_{i_{mn}}:=|y_{mn}|^2(\be_m-\be_n)(\be_m-\be_n)^\top +|y_{mn}|^2 (\be_{N_b+m}-\be_{N_b+n})(\be_{N_b+m}-\be_{N_b+n})^\top$.

The active power injected into bus $n$ can be decomposed into a dispatchable component $p_n^g$ and an inflexible component $p_n^d$ as $p_n=p_n^g-p_n^d$. The former captures the active power dispatch of a generator or a flexible load located at bus $n$. The latter captures the inelastic load to be served at bus $n$. To simplify the exposition, each bus is assumed to be hosting at most one dispatchable resource (generator or flexible load). The reactive power injected into bus $n$ is decomposed similarly as $q_n=q_n^g-q_n^d$. Let $\mcN_g\subseteq \mcN$ be the subset of buses hosting dispatchable power injections with cardinality $N_g$. Bus $n=1$ belongs to $\mcN_g$ and serves as the angle reference, so that $v_1^i=\bv^\top\be_{N_b+1}\be_{N_b+1}^\top\bv=0$. The remaining buses host non-flexible loads and constitute the subset $\mcN_\ell=\mcN\setminus \mcN_g$ with cardinality $N_l=N_b-N_g$. For simplicity, we will henceforth term the buses in $\mcN_g$ as \emph{generator buses}, and the ones in $\mcN_\ell$ as \emph{load buses}. Zero-injection (junction) buses belong to $\mcN_\ell$ and satisfy $p_n^g=p_n^d=q_n^g=q_n^d=0$.

Given the inflexible loads at all buses $\{p_n^d,q_n^d\}_{n\in\mcN}$, the OPF problem aims at optimally dispatching generators and flexible loads $\{p_n^g,q_n^g\}_{n\in\mcN_g}$ while meeting resource and network limits. The OPF can be formulated as the QCQP~\cite{Bose}, \cite{MSL15}
\begin{subequations}\label{eq:P1}
	\begin{align}
	\min\ &~\sum_{n\in\mcN_g} c_n^p p_n^g+c_n^q
	q_n^g\tag{P1}\label{eq:P1:cost}\\
	\mathrm{over}\ &~\bv\in\mathbb{R}^{2N_b},\{p_n^g,q_n^g\}_{n\in\mcN_g}\notag\\
	\mathrm{s.to}\ &~\bv^{\top}\bM_{p_n}\bv=p_n^g-p_n^d,&\forall~n\in\mcN_g\label{eq:P1:pg}\\
	&~\bv^{\top}\bM_{q_n}\bv=q_n^g-q_n^d,&\forall~n\in\mcN_g\label{eq:P1:qg}\\
	&~\bv^{\top}\bM_{p_n}\bv=-p_n^d,&\forall~n\in\mcN_\ell\label{eq:P1:pl}\\
	&~\bv^{\top}\bM_{q_n}\bv=-q_n^d,&\forall~n\in\mcN_\ell\label{eq:P1:ql}\\
	&~\underline{p}_n^g\leq \bv^{\top}\bM_{p_n}\bv+p_n^d\leq \bar{p}_n^g,&\forall~n\in\mcN_g\label{eq:P1:pglim}\\
	&~\underline{q}_n^g\leq \bv^{\top}\bM_{q_n}\bv+q_n^d\leq \bar{q}_n^g,&\forall~n\in\mcN_g\label{eq:P1:qglim}\\
	&~\underline{v}_n\leq \bv^{\top}\bM_{v_n}\bv\leq \bar{v}_n,&\forall~n\in\mcN\label{eq:P1:vlim}\\	
	&~\bv^{\top}\be_{N+1}\be_{N+1}^\top\bv=0\label{eq:P1:ref}\\	
	&~\bv^{\top}\bM_{i_{mn}}\bv\leq \bar{i}_{mn},&\forall~(m,n)\in\mcE\label{eq:P1:flim}
	\end{align} 
\end{subequations}
where $(c_n^p,c_n^q)$ are the coefficients for generation cost or the utility function for flexible load at bus $n$. Constraints \eqref{eq:P1:pg}--\eqref{eq:P1:ql} enforce the power flow equations at load and generator buses. Constraints \eqref{eq:P1:pglim}--\eqref{eq:P1:qglim} impose limits for generators and flexible loads. Constraints \eqref{eq:P1:vlim} confine squared voltage magnitudes within given ranges and \eqref{eq:P1:ref} identifies the reference bus. Finally, constraint \eqref{eq:P1:flim} limits squared current magnitudes according to line ratings.

Problem \eqref{eq:P1} is a parametric QCQP as it needs to be solved for different values of demands $\{p_n^d,q_n^d\}_{n\in\mcN}$; costs $\{c_n\}_{n\in\mcN_g}$; and generation capacities $\{\underline{p}_n^g,\bar{p}_n^g,\underline{q}_n^g,\bar{q}_n^g\}$. Voltage limits $\{\underline{v}_n,\bar{v}_n\}$ for $n\in\mcN$ and current limits $\{\bar{i}_{(m,n)}\}$ for $(m,n)\in\mcE$ may also be changing due to normal and emergency ratings. To keep the exposition uncluttered, we henceforth fix all but the inelastic demands to known values. In other words, we are interested in solving \eqref{eq:P1} over different values of the parameter vector $\btheta:=\{p_n^d,q_n^d\}_{n\in\mcN}\in\mathbb{R}^{2L}$. 

The optimization variables of \eqref{eq:P1} consist of all nodal voltages $\bv$ and the (re)active power schedules for generators. If vector $\bx_g$ collects generator schedules $\{p_n^g,q_n^g\}_{n\in\mcN_g}$, then vector $\bx_{\btheta}^\top:=[\bv^\top~\bx_g^\top]$ denotes the minimizer of \eqref{eq:P1} for the specific parameter vector $\btheta$. Aiming for the complete $\bx_{\btheta}$ is apparently an  over-parameterization of the problem, adopted only to ease the formulation in \eqref{eq:P1}. What the system operator actually needs to know in practice is only the voltage magnitude and active power schedule for each generator (modulo the reference generator for which we set the voltage magnitude and angle). In light of this and to reduce the DNN output dimension, the DNN is trained to predict the PV setpoints for generators. Given the predicted quantities and knowing the values for inflexible loads from $\btheta$, the remaining quantities can be readily computed using a power flow solver anyway.

Given a set of (locally) optimal primal/dual solutions for~\eqref{eq:P1}, we next analyze the sensitivity of the parametric AC-OPF. Sensitivities are computed for the complete $\bx_{\btheta}$, from which the Jacobian $\bJ_{\btheta}$ needed in \eqref{eq:si-dnn} can be readily obtained.


\section{Sensitivity Analysis for QCQP-based OPF}\label{sec:sensitivity}
To analyze the sensitivity of \eqref{eq:P1} with respect to $\btheta$, let us first express the vectors of complex power injections across all buses as
\begin{equation}\label{eq:AB}
\left[\begin{array}{c}
\bp\\
\bq
\end{array}\right]=\bA\bx_g+\bB\btheta
\end{equation}
where $\bx_g$ stacks the generator (re)active power injections $\{p_n^g,q_n^g\}_{n\in\mcN_g}$ and $(\bA,\bB)$ are matrices assigning generators and loads to buses. We then reformulate \eqref{eq:P1} as
\begin{subequations}\label{eq:P2}
	\begin{align}
	\min_{\bv,\bx_g}\ &~\ba_0^\top\bx_g\label{eq:P2:cost}\\
	\mathrm{s.to}\ &~\bv^{\top}\bL_{\ell}\bv=\ba_{\ell}^\top \bx_g+\bb_{\ell}^\top \btheta,&~\ell=1:L:~~\lambda_{\ell}\label{eq:P2:g}\\
	&~\bv^{\top}\bM_m\bv\leq\bd_m^\top \btheta+f_m,&~m=1:M:~~\mu_m\label{eq:P2:m}
	\end{align} 
\end{subequations}
where $\ba_0$ collects the generation cost coefficients (cf.~\eqref{eq:P1}); the first $L=2N_b+1$ constraints in \eqref{eq:P2:g} correspond to the power flow equations \eqref{eq:P1:pg}--\eqref{eq:P1:ql} and the angle reference constraint \eqref{eq:P1:ref}, while constraints \eqref{eq:P2:m} correspond to the $M=4N_g+2N_b+E$ inequality constraints of \eqref{eq:P1}. Matrices $(\bL_{\ell},\bM_m)$ are drawn from the $\bM$ matrices appearing in the quadratic forms of \eqref{eq:P1}. Vectors $(\ba_{\ell},\bb_{\ell})$ correspond to rows of matrices $(\bA,\bB)$ in \eqref{eq:AB} for $\ell\leq2N_b$; and $\bzero$ for $\ell=2N_b+1$. Vectors $\bd_m$ are indicator (canonical) vectors and constants $f_m$ relate to generation, voltage, and line limits. 

Aiming at computing the sensitivity of a minimizer $\bx_{\btheta}^\top=[\bv^\top~\bx_g^\top]$ of~\eqref{eq:P2} with respect to $\btheta$, we explored the related literature. There has indeed been significant interest in computing the sensitivities of OPF minimizers with respect to load~\cite{Almeida94parametric}, \cite{Ajjarapu95OCPF}, \cite{Almeida2000varload}. However, the primary motivation for these works was to efficiently compute minimizers and look into binding constraints for a \emph{given trajectory} of load variations. Hence the related OPF was parameterized using a scalar conveniently varied over a range of interest. Seeking to compute the minimizer sensitivities with respect to the vector $\btheta$ in a relatively general setting, we explored beyond the power systems literature. Fortunately, there exists a rich corpus of work on perturbation analysis of continuous optimization problems with applications in operation research, economics, mechanics, and optimal control~\cite{shapiro2013perturbation}. The first approaches applied the implicit function theorem to the related first-order optimality conditions~\cite{fiacco1976sensitivity}. Thereon, many developments have been made towards relaxing the assumptions of initial works, and expanding the scope to conic programs~\cite{Conejo06}, \cite{Castillo06}, \cite{shapiro2013perturbation}, \cite{agrawal2020differentiating}. For several recent applications however, the early approaches of~\cite{fiacco1976sensitivity} are well suited due to their simplicity; see for example~\cite{AmosKotler17}. Building upon~\cite{fiacco1976sensitivity}, we next compute the sensitivities required for SI-DNN in Section~\ref{subsec:perturb}; and relax some of the needed assumptions in Section~\ref{subsec:existence}.

\subsection{Perturbing Optimal Primal/Dual Solutions}\label{subsec:perturb}
Towards instantiating \eqref{eq:Ptheta} with~\eqref{eq:P2} and to reduce notational clutter, let us use symbols $(\bx,\blambda,\bmu)$ to denote the \emph{optimal} primal/dual variables $(\bx_{\btheta},\blambda_{\btheta},\bmu_{\btheta})$ of \eqref{eq:P2} for a particular $\btheta$. Under mild technical assumptions, a local primal/dual point for \eqref{eq:P2} satisfies the first-order optimality conditions~\cite{Conejo06}. The goal of sensitivity analysis is to find infinitesimal changes $(\d{\bx},\d\blambda,\d\bmu)$, so that the perturbed point $(\bx+\d\bx, \blambda+\d\blambda,\bmu+\d\bmu)$ still satisfies the first-order optimality conditions when the input parameters change from $\btheta$ to $\btheta+\d\btheta$~\cite{fiacco1976sensitivity}. To this end, we next review the optimality conditions and then differentiate them to compute the sought sensitivities.

The Lagrangian function of \eqref{eq:P2} is defined as
\begin{align*}
    \mcL(\bx, \blambda,\bmu;\btheta):=\ba_0^\top\bx_g&+\sum_{\ell=1}^L\lambda_\ell\left(\bv^{\top}\bL_{\ell}\bv-\ba_{\ell}^\top \bx_g-\bb_{\ell}^\top \btheta\right)\\
    &+\sum_{m=1}^M\mu_m\left(\bv^{\top}\bM_m\bv-\bd_m^\top \btheta-f_m\right).
\end{align*}
With $\bx:=\{\bv,\bx_g\}$, Lagrangian optimality $\nabla_{\bx}\mcL=\bzero$ gives
\begin{subequations}\label{eq:kkt1}
    \begin{align}
     &\underbrace{\left(\sum_{\ell=1}^L \lambda_\ell \bL_\ell +\sum_{m=1}^M \mu_m\bM_m\right)}_{:=\bZ}\bv=\bzero.\label{eq:kkt1:a}\\
    &\ba_0=\sum_{\ell=1}^L\lambda_{\ell} \ba_{\ell}\label{eq:kkt1:b}
    \end{align}
\end{subequations}
In addition to Lagrangian optimality, first-order optimality conditions include primal feasibility [cf.~\eqref{eq:P2:g}--\eqref{eq:P2:m}], as well as complementary slackness and dual feasibility for all $m$:
\begin{subequations}\label{eq:kkt2}
    \begin{align}
    &\mu_m\underbrace{\left(\bv^{\top}\bM_m\bv-\bd_m^\top \btheta-f_m\right)}_{:=g_m}=0\label{eq:kkt2:a}\\
    &\mu_m \geq 0.\label{eq:kkt2:b}
    \end{align}
\end{subequations}
From the aforementioned optimality conditions, let us focus on those that take the form of equalities, namely \eqref{eq:kkt1:a}--\eqref{eq:kkt1:b}, \eqref{eq:P2:g}, and \eqref{eq:kkt2:a}. For these conditions, we will compute their total differentials. From the first three, we obtain
\begin{subequations}\label{eq:diff1}
    \begin{align}
    \bZ\d\bv+\bL_{\lambda}\d \blambda+\bM_{\mu}\d \bmu&=\bzero\label{eq:diff2:a}\\
    \bA^\top\d\blambda&=\bzero\label{eq:diff2:b}\\
    2\bL_{\lambda}^\top\d\bv-\bA\d \bx_g -\bB\d\btheta&=\bzero\label{eq:diff2:c}
    \end{align}
\end{subequations}
where $\bL_{\lambda}:=\sum_{\ell=1}^L\bL_\ell\bv\be_\ell^\top$; and $\bM_{\mu}:=\sum_{m=1}^M\bM_m\bv\be_m^\top$.

For~\eqref{eq:kkt2:a}, the total differential is 
\begin{equation}\label{eq:cs}
    g_m\d \mu_m + \mu_m \d g_m=0
\end{equation}
where $\d g_m := \left(\nabla_\bv g_m\right)^\top\d \bv + \left(\nabla_{\btheta} g_m\right)^\top\d \btheta$ for all $m$. We identify three cases: 
\renewcommand{\labelenumi}{\emph{c\arabic{enumi})}}
\begin{enumerate}
    \item If $\mu_m=0$ and $g_m<0$, then \eqref{eq:cs} implies $\d \mu_m=0$. It follows that: \emph{i)} $\mu_m+\d \mu_m=0$; \emph{ii)} $(\mu_m+\d \mu_m)(g_m+\d g_m)=0$; and \emph{iii)} $g_m+\d g_m<0$ for any small $\d g_m$. In conclusion, condition~\eqref{eq:cs} ensures that the perturbed point satisfies conditions for optimality, including the inequalities from primal/dual feasibility.
    \item If $\mu_m>0$ and $g_m=0$, then \eqref{eq:cs} gives $\d g_m=0$. It also follows that: \emph{i)} $g_m+\d g_m=0$; \emph{ii)} $(\mu_m+\d \mu_m)(g_m+\d g_m)=0$; and \emph{iii)} $\mu_m+\d \mu_m>0$ for any small $\d \mu_m$. As in case \emph{c1)}, condition~\eqref{eq:cs} ensures that the perturbed point satisfies all conditions for optimality.
    \item If $\mu_m=g_m=0$, then \eqref{eq:cs} is inconclusive on $\d g_m$ and $\d \mu_m$. In this \emph{degenerate} case, for the perturbed point to remain optimal, we need to explicitly impose: \emph{i)} $\d g_m\leq 0$; \emph{ii)} $\d \mu_m\geq0$; and \emph{iii)} $\d g_m \d \mu_m=0$. Even though the three latter constraints can be handled by the sensitivity analysis of~\cite{Conejo06}, \cite{Castillo06}, they considerably complicate the treatment. Moreover, such degeneracy is seldom encountered numerically. We henceforth rely on the so called \emph{strict complementarity} assumption, which ignores case \emph{c3)}~\cite{fiacco1976sensitivity}.
\end{enumerate}

\begin{assumption}\label{as:scs}
Given a tuple of optimal primal/dual variables $(\bx,\blambda,\bmu)$, constraint $g_m(\bx;\btheta)=0$ if and only if $\mu_m>0$.
\end{assumption}

Two observations are in order. First, the analysis under \emph{c1)-c2)} reveals that although we perturbed only the equality conditions for optimality, the obtained perturbed point satisfies the inequality conditions for optimality as well. Therefore, under Assumption~\ref{as:scs}, the point $(\bx+\d\bx, \blambda+\d\blambda,\bmu+\d\bmu)$ satisfying the perturbed optimality conditions is (locally) optimal for \eqref{eq:P2}, when solved for ${\btheta+\d\btheta}$. Second, despite Assumption~\ref{as:scs}, if a degenerate instance of \eqref{eq:P2} does occur for some $\btheta_s$ in the training dataset, the particular pair $(\btheta_s,\bx_{\btheta_s})$ can still be used to train the SI-DNN, yet without the additional sensitivity information. In other words, degenerate instances can contribute only to the first fitting term of \eqref{eq:si-dnn}.

Applying \eqref{eq:cs} for all $m$, the total derivatives for \eqref{eq:kkt2:a} can be compactly expressed as
\begin{equation}
    \mcD(\bg) \d \bmu +2\mcD(\bmu)\bM_{\mu}^\top 
        \d\bv -\bD^\top \d\btheta=\bzero,\label{eq:diff:cs}
\end{equation}
where $\bg:=\{g_m\}_{m=1}^M$ stacks the inequality constraint values, and matrix $\bD:=\sum_{m=1}^M\mu_m\bd_m\be_m^\top$. Operator $\mcD(\bx)$ returns a diagonal matrix with vector $\bx$ on its main diagonal. Conditions \eqref{eq:diff1} and \eqref{eq:diff:cs} can be collected in matrix-vector form as
\begin{equation}\label{eq:diff}
    \underbrace{\left[\begin{array}{cccc}
        \bZ & \bzero & \bL_{\lambda} & \bM_{\mu} \\
        \bzero & \bzero & \bA^\top &\bzero\\
        2\bL_{\lambda}^\top & -\bA & \bzero & \bzero \\
        2\mcD(\bmu)\bM_\mu^\top & \bzero & \bzero & \mcD(\bg)
    \end{array}\right]}_{:=\bS}
       \left[\begin{array}{c}
        \d \bv\\
        \d \bx_g\\
        \d \blambda\\
        \d \bmu
    \end{array}\right] {=}
        \underbrace{\left[\begin{array}{c}
        \bzero\\
        \bzero\\
        \bB\\
        \bD^\top    
        \end{array}\right]}_{:=\bU}\d \btheta
\end{equation}

To compute the sensitivities of primal and dual variables with respect to the $p$-th entry $\theta_p$ of $\btheta$, we need to solve the previous system of $2(N_b+N_g)+L+M$ linear equations for $\d\btheta=\be_p$. The size of the system can be reduced by dropping the numerous inactive inequality constraints of \eqref{eq:P2} for which $\mu_m=0$ and $g_m<0$, and thus, $\d \mu_m=0$ as discussed earlier under case \emph{c1)}. Notably, if matrix $\bS$ is invertible, the aforementioned sensitivities can all be found at once using the respective blocks of $\bS^{-1}\bU\be_p$. We next address two relevant questions: \emph{q1) When is $\bS$ invertible?}; and \emph{q2) What are the implications of a singular $\bS$?}

\subsection{Existence of Primal Sensitivities}\label{subsec:existence}
To address \emph{q1)} for an arbitrary \eqref{eq:Ptheta}, the existing literature identifies some assumptions on $(\bx,\blambda,\bmu;\btheta)$. We first review these assumptions, and then assess if they are reasonable for the OPF task at hand. Given an optimal primal $\bx$ for some $\btheta$, let $\mcA(\bx)$ denote the subset of inequality constraints of $\bg(\bx;\btheta)\leq\bzero$ that are \emph{active or binding}, that is $\mcA(\bx):=\{m:g_m(\bx;\btheta)=0\}$. A primal solution $\bx$ is termed \emph{regular} if the next assumption holds.

\begin{assumption}\label{as:regular}
The vectors $\{\nabla_{\bx}h_\ell\}_{\forall \ell}$ and $\{\nabla_{\bx}g_m\}_{m\in\mcA(\bx)}$ are linearly independent.
\end{assumption}

For the OPF in~\eqref{eq:P2}, the functions $h_\ell$ and $g_m$ correspond to the (in)equality constraints \eqref{eq:P2:g}--\eqref{eq:P2:m} written in the standard form as in~\eqref{eq:Ptheta}. Assumption~\ref{as:regular} is often referred to as linearly independent constraint qualification~(LICQ). If a (locally) optimal $\bx$ satisfies the LICQ, the corresponding optimal dual variables $(\blambda,\bmu)$ are known to be unique~\cite{Be99}. In addition to satisfying first-order optimality conditions, a sufficient condition for $(\bx,\blambda,\bmu;\btheta)$ to be (locally) optimal is often provided by the following second-order optimality condition.

\begin{assumption}\label{as:secorder}
For a subspace orthogonal to the subspace spanned by the gradients of active constraints 
\[\mcZ:=\left\{\bz: \bz^\top\nabla_{\bx}h_\ell=0~\forall~\ell,\bz^\top\nabla_{\bx}g_m=0~\forall~m\in\mcA(\bx)\right\}\]
it holds that $\bz^\top\nabla_{\bx\bx}^2\mcL\bz>0$ for all $z\in\mcZ\setminus\{\bzero\}$.
\end{assumption}

Under the strict complementarity, regularity, and second-order optimality conditions, matrix $\bS$ is guaranteed to be invertible; see Theorem~2.1 and Corollary~2.1 of \cite{fiacco1976sensitivity}.

\begin{lemma}[\cite{fiacco1976sensitivity}]\label{le:Sinv}
If Assumptions~\ref{as:scs}--\ref{as:secorder} hold, matrix $\bS^{-1}$ exists.
\end{lemma}

Lemma~\ref{le:Sinv} implies that under the stated assumptions, the optimal primal and dual variables of \eqref{eq:Ptheta} vary smoothly with changes in parameter $\btheta$, and the associated sensitivities can be found via~\eqref{eq:diff}. Prior works that compute sensitivities of optimal primal and dual variables for scalar-parameterized OPF instances rely on the non-singularity of $\bS$~\cite{Almeida94parametric}, \cite{Ajjarapu95OCPF}.

While Assumptions~\ref{as:scs}--\ref{as:secorder} seem to be standard in the optimization literature, our recent work on the optimal dispatch of inverters in distribution grids demonstrated analytically and numerically that LICQ (Assumption~\ref{as:regular}) is violated frequently~\cite{SGKCB2020}. Instances violating LICQ can be conceived for the AC-OPF in~\eqref{eq:P2} too~\cite{Almeida16ill}, \cite{Dorfler18LICQ}. To bring up one such example, consider a power system where a load bus $m$ is connected to the rest of the system through another bus $n$ via a single transmission line $(m,n)$. As bus $m$ is a load bus, it contributes two equality constraints \eqref{eq:P1:pl} and \eqref{eq:P1:ql}. It can be shown that if any of the three following scenarios occurs, LICQ fails: \emph{i)} the voltage limits in \eqref{eq:P1:vlim} become binding (above or below) for both buses $m$ and $n$; \emph{ii)} line $(m,n)$ becomes congested [cf.~\eqref{eq:P1:flim}] and a voltage limit at bus $m$ becomes binding; or \emph{iii)} line $(m,n)$ becomes congested and a voltage limit at bus $n$ becomes binding. Further detailed examples for AC-OPF instances violating LICQ can be found in~\cite{Almeida16ill}, \cite{Dorfler18LICQ}. Attempting to circumvent LICQ violation via problem reformulations may be futile as their occurrences depend on $\btheta$, and are thus hard to analyze. Under certain assumptions on OPF instances and load variations, LICQ occurrences can be shown to have zero measure~\cite{Dorfler18LICQ}. If the required assumptions are not met, resorting to Fritz-John rather than the KKT conditions for sensitivity analysis has been proposed~\cite{Almeida16ill}. However, before tackling the singularity of $\bS$ due to LICQ violation, we must answer question \emph{q2)}. 

The implications of a singular $\bS$ have previously been investigated in \cite{Conejo06} and \cite{Castillo06}: When LICQ is violated despite strict complementarity, the sensitivities of some primal/dual variables may still exist with respect to a $\theta_p$. In detail, consider the set $\Gamma:=\{\bgamma\in\mathbb{R}^{2(N_b+N_g)+L+M}:\bS\bgamma=\bU\be_p\}$, which is the solution set of \eqref{eq:diff}. If the $n$-th entry of $\bgamma$ remains constant for all $\bgamma\in\Gamma$, the sensitivity of the $n$-th entry of $[\bx^\top~\blambda^\top~\bmu^\top]^\top$ with respect to $\theta_p$ does exist; see~\cite{Conejo06} and \cite{Castillo06} for physical interpretation and illustrative examples. While a subset of optimal primal/dual variables may be differentiable under LICQ violation, explicitly identifying the differentiable quantities requires instance-based numerical evaluation in ~\cite{Conejo06} and \cite{Castillo06}. Since for training an SI-DNN, we are interested only in the sensitivities $\nabla_{\btheta}\bx$, we need to ensure that all solutions $\bgamma\in\Gamma$ share the same first $N$ entries. This is equivalent to saying that the first entries of $\bn$ are zero for all $\bn\in\nullspace(\bS)$. The equivalence stems from the fact that if $\bS\bbgamma=\bu$ for a $\bbgamma$, any other solution to $\bS\bgamma=\bu$ takes the form $\bgamma=\bbgamma+\bn$ for some $\bn\in\nullspace(\bS)$. The next claim provides sufficient conditions for the first $N$ entries of $\bn$ to be zero. 


\begin{theorem}\label{th:th1}
If Assumptions~\ref{as:scs} and~\ref{as:secorder} hold, then $n_i=0$ for $i=1,\dots,2(N_b+N_g)$ for all $\bn\in\nullspace(\bS)$.
\end{theorem}

\begin{IEEEproof}
The claim holds trivially for $\bn=\bzero$. The proof for non-zero $\bn$ builds on Assumption~\ref{as:secorder}, and thus the terms involved in these assumptions are computed for~\eqref{eq:P2} first.
	\begin{align*}
	\nabla_{\bx\bx}^2&\mcL:=\begin{bmatrix}\bZ&\bzero\\ \bzero&\bzero\end{bmatrix},\nabla_{\bx}h_\ell:=\begin{bmatrix}2\bL_{\ell}\bv\\-\ba_{\ell}\end{bmatrix},
	\nabla_{\bx}g_m:=\begin{bmatrix}2\bM_m\bv\\\bzero\end{bmatrix}.
	\end{align*} 
Given vector $\bn\neq\bzero$ such that $\bS\bn=\bzero$, partition $\bn$ conformably to $[\bv^\top~\bx_g^\top~\blambda^\top~\bmu^\top]^\top$ as  $\bn:=[\bn_v^\top~\bn_{x_g}^\top~\bn_{\lambda}^\top~\bn_{\mu}^\top]^\top$. By the definition of $\bS$ in~\eqref{eq:diff}, expanding $\bS\bn=\bzero$ yields
\begin{subequations}\label{eq:expand}
    \begin{align}
    \bZ\bn_v + \bL_{\lambda}\bn_{\lambda}+ \bM_{\mu}\bn_{\mu}&=\bzero\label{eq:expand:a}\\
    \bA^\top\bn_{\lambda}&=\bzero\label{eq:expand:b}\\
    2\bL_{\lambda}^\top\bn_v-\bA\bn_{x_g}&=\bzero\label{eq:expand:c}\\
    2\mcD(\bmu)\bM_{\mu}^\top\bn_{v}+ \mcD(\bg)\bn_{\mu}&=\bzero.\label{eq:expand:d}
    \end{align}
\end{subequations}
Assumption~\ref{as:scs} dictates that the first and second term in~\eqref{eq:expand:d} have complementary sparsity, so that $\mcD(\bmu)\bM_{\mu}^\top\bn_{v}=\bzero$ and $\mcD(\bg)\bn_{\mu}=\bzero$. Recalling the definition $\bM_{\mu}:=\sum_{m=1}^M\bM_m\bv\be_m^\top$, equation $\mcD(\bmu)\bM_{\mu}^\top\bn_{v}=\bzero$ indeed implies $[\bn_{v}^\top~\bn_{x_g}^\top]\nabla_{\bx}g_m=0$ for all $m\in\mcA(\bx)$, and together with \eqref{eq:expand:c} ensures
\begin{equation}\label{eq:ninZ}
    \begin{bmatrix}\bn_{v}\\\bn_{x_g}\end{bmatrix}\in\mcZ.
\end{equation}
Pre-multiplying \eqref{eq:expand:a}--\eqref{eq:expand:b} by $2\bn_{v}^\top$ and $\bn_{x_g}^\top$ and subtracting the two resulting equations yields
\begin{align}\label{eq:sum_ab}
    2\bn_{v}^\top\bZ\bn_{v} + 2\bn_{v}^\top\bL_{\lambda}\bn_{\lambda}-\bn_{x_g}^\top\bA^\top\bn_{\lambda}\notag\\+ [\bn_v^\top~\bn_x^\top]\begin{bmatrix}
    2\bM_{\mu}\\\bzero
    \end{bmatrix}\bn_{\mu}=0
\end{align}
where the second and third term on the left-hand side (LHS) sum up to zero per~\eqref{eq:expand:c}. Since $\mcD(\bg)\bn_{\mu}=\bzero$, if $g_m\neq0$, then the $m$-th entry $n_{\mu,m}$ of $\bn_{\mu}$ should be zero. In other words, we get that $n_{\mu,m}=0$ for all $m\notin\mcA(\bx)$, and thus, $2[\bM_{\mu}^\top~\bzero]^\top\bn_{\mu}=\sum_{m\in\mcA(\bx)}\nabla_{\bx}g_mn_{\mu,m}$. Substituting the latter into \eqref{eq:sum_ab} gives
\begin{align*}
    2\bn_{v}^\top\bZ\bn_{v} + \sum_{m\in\mcA(\bx)}[\bn_{v}^\top~\bn_{x_g}^\top]^\top\nabla_{\bx}g_mn_{\mu,m}=0.
\end{align*}
The second term on the LHS equals zero due to~\eqref{eq:ninZ}. Therefore $\bn_{v}^\top\bZ\bn_{v}=0$, thus implying
$$\begin{bmatrix}\bn_{v}^\top &\bn_{x_g}^\top\end{bmatrix}\nabla_{xx}^2\mcL\begin{bmatrix}\bn_{v}\\\bn_{x_g}\end{bmatrix}=0 $$ which contradicts Assumption~\ref{as:secorder}, unless $\bn_{v}=\bzero$ and $\bn_{x_g}=\bzero$.
\end{IEEEproof}

Thanks to Theorem~\ref{th:th1}, we can proceed with computing $\nabla_{\btheta}\bx$ by solving~\eqref{eq:diff} even if $\bS$ is singular. In other words, Theorem~\ref{th:th1} allows us to compute $\nabla_{\btheta}\bx$ even if the LICQ (Assumption~\ref{as:regular}) fails. If $\bS^\dagger$ is the pseudo-inverse of $\bS$, the Jacobian matrix $\bJ_{\btheta}=\nabla_{\btheta}\bx$ can be computed as the top $2(N_b+N_g)$ rows of $-\bS^\dagger\bU$. 

The previous analysis has tacitly presumed the system $\bS\bgamma=\bu$ has at least one solution for all $\bu\in\range(\bU)$. The numerical tests of Section~\ref{sec:tests} demonstrate that for the AC-OPF in~\eqref{eq:P2}, the system $\bS\bgamma=\bu$ features a solution indeed.

As discussed earlier, we focus on training an SI-DNN for predicting generator voltage magnitudes and active power setpoints. Having solved~\eqref{eq:diff} and found the sensitivity of~$\bx$ with respect to $\btheta$, the sensitivity of active power generation can be obtained readily using the corresponding entries of $\bx_g$. The sensitivity of voltage magnitudes can be derived from the sensitivities of the real and imaginary components of voltages with respect to $\btheta$. Precisely, the voltage magnitude at bus~$n$ is given by $v_n=\sqrt{(v_n^r)^2+(v_n^i)^2}$ and its sensitivity with respect to $\theta_\ell$ can be found through the chain rule
$$\frac{\partial v_n}{\partial \theta_\ell}=\frac{1}{v_n}\left(v_n^r\frac{\partial v_n^r}{\partial \theta_\ell}+v_n^i\frac{\partial v_n^i}{\partial \theta_\ell}\right).$$
Evaluating the above completes the requirements of sensitivities for augmenting the SI-DNN training set.

From \eqref{eq:diff}, the required sensitivities obviously depend on values of optimal primal/dual variables of \eqref{eq:P2}. Problem \eqref{eq:P2} is a non-convex quadratic program and existing solvers may converge to a local rather than a global solution. Albeit the previous sensitivity analysis is valid even for local solutions, the performance of the trained DNN will be apparently suboptimal. To train an SI-DNN to predict global OPF solutions, we next extend the analysis to the SDP relaxation of \eqref{eq:P2}.

\section{SDP Relaxation of the AC-OPF}\label{sec:SDP}
In pursuit of globally optimal AC-OPF schedules, the non-convex QCQP of \eqref{eq:P2} can be relaxed to the SDP~\cite{Bai08}
\begin{subequations}\label{eq:SDP}
	\begin{align}
	\min_{\bx_g,\bV\succeq 0}\ &~\ba_0^\top\bx_g\label{eq:SDP:cost}\\
	\mathrm{s.to}\ &~\trace(\bL_{\ell}\bV)=\ba_{\ell}^\top \bx_g+\bb_{\ell}^\top \btheta,&~\forall \ell:~\lambda_{\ell}\label{eq:SDP:g}\\
	&~\trace(\bM_m\bV)\leq\bd_m^\top \btheta+f_m,&~\forall m:~~\mu_m.\label{eq:SDP:m}
	\end{align} 
\end{subequations}
Problem~\eqref{eq:SDP} is equivalent to \eqref{eq:P2} if matrix $\bV$ is rank-1 at optimality, in which case the SDP relaxation is deemed as \emph{exact}. The relaxation turns out to be exact for several power networks and practical loading conditions; see~\cite{Low14} for a review of related analyses. When the relaxation is exact, the $\bV$ minimizer of \eqref{eq:SDP} can be expressed as $\bV=\bv\bv^\top$. 

We briefly review how the solution $(\bx_g,\bv;\blambda,\bmu)$ obtained from the SDP formulation of \eqref{eq:SDP} satisfies the first-order optimality conditions for the non-convex QCQP in \eqref{eq:P2} as well. To this end, it is not hard to derive the dual program of \eqref{eq:SDP}:
\begin{subequations}\label{eq:SDPd}
	\begin{align}
	\max_{\blambda,\bmu}&\ -\sum_{\ell=1}^L\lambda_{\ell}\bb_{\ell}^\top\btheta-\sum_{m=1}^M\mu_m(\bd_m^\top\btheta+f_m)\label{eq:SDPd:cost}\\
	\mathrm{s.to}\ &~\ba_0=\sum_{\ell=1}^L\lambda_{\ell}\ba_{\ell}\label{eq:SDPd:con:1}\\ &~\bZ:=\sum_{\ell=1}^L\lambda_{\ell}\bL_{\ell}+\sum_{m=1}^M\mu_m\bM_m\succeq 0\label{eq:SDPd:con:2}\\
	&~\bmu\geq \bzero.\label{eq:SDPd:con:3}
	\end{align} 
\end{subequations}
The optimality conditions for the SDP primal-dual pair \eqref{eq:SDP}--\eqref{eq:SDPd} then include:
\renewcommand{\labelenumi}{\emph{\roman{enumi})}}
\begin{enumerate}
    \item Primal feasibility~\eqref{eq:SDP:g}--\eqref{eq:SDP:m} implies~\eqref{eq:P2:g}--\eqref{eq:P2:m}.
    \item Dual feasibility~\eqref{eq:SDPd:con:3} applies to~\eqref{eq:P2} as well.
    \item Complementary slackness for~\eqref{eq:SDP:m} applies to~\eqref{eq:P2:m}.
    \item Complementary slackness for~\eqref{eq:SDPd:con:2} gives $\trace(\bV\bZ)=0$ or $\bZ\bv=\bzero$, which along with~\eqref{eq:SDPd:con:1} yield the Lagrangian optimality conditions for the QCQP shown in \eqref{eq:kkt1}.
\end{enumerate}
Therefore $(\bx_g,\bv;\blambda,\bmu)$ is a stationary point for the QCQP in \eqref{eq:P2}. Because it further attains the optimal cost for the relaxed problem in \eqref{eq:SDP}, it is in fact the globally optimal for \eqref{eq:P2}. 

To recapitulate, we have used the non-convex QCQP formulation of the AC-OPF to derive the sensitivity formulae of \eqref{eq:diff}. This is advantageous as the QCQP features differentiable objective and constraint functions. The obtained sensitivity formulae can be evaluated at the AC-OPF solution provided by any nonlinear programming solver, although such solution may be only locally optimal. To compute a globally optimum AC-OPF solution, we propose using \eqref{eq:SDP} instead. If the SDP relaxation is exact, the obtained solution is globally optimal, while the sensitivity formulae derived from QCQP can still be used. Our suggested workflow avoids computing the sensitivities of the SDP formulation for the AC-OPF: Even though differentiating through convex cone constraints is possible~\cite{agrawal2020differentiating}, it can be perplexing.

\begin{remark}\label{re:molzhan}
The aforesaid workflow runs the SDP-based solver to obtain an OPF solution, but computes its sensitivities using the convenient formulae associated with the QCQP-OPF. This is to ensure global optimality if the SDP relaxation is exact. An alternative way to check global optimality is to follow the workflow of~\cite{Molzhan14}: Obtain an OPF solution via a mature OPF solver (e.g., QCQP or MATPOWER), and use the optimality conditions of the SDP-based OPF to check whether the obtained QCQP- or MATPOWER-based solution is globally optimal. Nevertheless, this optimality check relies on sufficient conditions. As a result, if the QCQP- or MATPOWER-based solution does not pass the global optimality test (that is indeed the case for the IEEE 300-bus system~\cite{Molzhan14}), one may still have to run the SDP-based OPF solver in pursuit of a better solution or a global optimality guarantee. 
\end{remark}

\section{Numerical Tests}\label{sec:tests}
The novel SI-DNN approach was evaluated using the IEEE 39-bus, the IEEE 118-bus, and the Illinois 200-bus system. Datasets were generated using either the nonlinear OPF MATPOWER or the globally optimal SDP-based solver. 

\subsection{DNN Architecture and Training}
To ease the implementation and without loss of generality, we assumed that buses hosting generators do not host loads, i.e., $p_n^d=q_n^d=0$ for all $n\in\mcN_g$. As discussed at the end of Section~\ref{sec:QCQP}, the DNN input~$\btheta$ consists of the $2(N_b-N_g)$ (re)active power demands at load buses. The DNN output is the setpoints for active power and voltage magnitude~$(p_n^g,v_n)$ at all generators $n\in\mcN_g$ excluding $p_1^g$ for the slack bus. We collect these output quantities in~$\cbx_{\btheta}$, a subvector of $\bx_{\btheta}$. 

Both for P-DNN and SI-DNN, we chose a feed-forward fully-connected architecture. For the number of hidden layers being $K$, denote the number of neurons in layer $k$ by $u_k$, with input dimension $u_0=2(N_b-N_g)$ and output dimension $u_{K+1}=2N_g-1$. To explicitly constrain DNN outputs as per~\eqref{eq:P1:pglim} and \eqref{eq:P1:vlim}, the output layer uses $\tanh$ as its activation function, while all other layers use $\mathrm{ReLU}$. For DNN training and evaluation, labels $\cbx_{\btheta}$ were suitably scaled within $[-1,1]$. 

We built all DNNs using the \texttt{TensorFlow 2.0} python platform alongside \texttt{Keras} libraries. Training an SI-DNN deviates from the default routine as gradient updates are implemented separately; see Appendix~\ref{sec:AppA} for key differences. For DNN training, at every weight-update step, the gradients computed via the procedure in the appendix are passed to the Adam optimizer. For all tests, optimizer Adam was used with an exponential decay reducing the rate to $85\%$ every $250$ epochs. The initial learning rate will be reported later. DNNs were compiled using Jupyter notebook on a 2.7 GHz Intel Core i5 computer with 8 GB RAM.


\subsection{Learning Locally Optimal OPF Solutions}\label{subsec:localtests}
We first trained DNNs towards predicting MATPOWER AC-OPF minimizers. We contrasted SI-DNN with P-DNN in terms of the MSE and the related training times. With the primary goal of improving sample efficiency, the numerical tests emphasize on performance evaluation for relatively small training datasets. Nevertheless, to gain insight on the effect of the training dataset size, we first present tests using larger training datasets. 

\subsubsection{Tests on IEEE 39-bus system with large training datasets}
The network parameters and nominal loads for the IEEE 39-bus system were fetched from MATPOWER {\verb|casefile|}~\cite{MATPOWER}. The 39-bus system hosts $N_g=10$ generators. The benchmark system has loads on two of the generator buses. Removing these, there are $29$ load buses. To build a dataset for DNN testing and training, a set of 12,000 random $\btheta\in\mathbb{R}^{2\cdot 29}$ was sampled. The corresponding $\cbx_{\theta}$'s were obtained via MATPOWER. The dataset thus obtained was partitioned into training, cross-validation, and testing sets of sizes 10,000; 1,000; and 1,000, respectively. If infeasibility is encountered for some $\btheta$'s, such instances were omitted from the dataset. To represent various demand levels, we sampled the 12,000 random $\btheta$'s by scaling the benchmark demands entry-wise by a scalar drawn independently and uniformly within $[0.8, 1.2]$. For the aforementioned sampling, all 12,000 OPF instances were feasible. Since the generator cost functions are identical in the benchmark system, a uniform active power cost was used for all generators. The default OPF formulation of MATPOWER deviates from the QCQP in~\eqref{eq:P1:cost}. These differences introduce some nuances in building the linear system of~\eqref{eq:diff} for computing sensitivities; see Appendix~\ref{sec:AppB} for details. Having built the aforementioned dataset $\{(\btheta_s,\bJ_{\btheta_s},\cbx_{\btheta_s})\}_{s=1}^{12000}$, the architectures for SI-DNN and P-DNN were determined next. Based on preliminary tests, identical architectures were chosen for P-DNN and SI-DNN with $K=4$ hidden layers with $u_k=256$ neurons for ${k=1,\dots,4}$. Preliminary tests showed negligible effect on P-DNN performance if the number of layers is reduced to three. Nevertheless, the architecture for the two DNNs was kept identical to ensure equal expressibility.

\begin{table}[t]
	\caption{Average Test MSE [$\times~10^{-3}$] for predicting MATPOWER solution on IEEE 39-bus system, and weighting factor $\rho$}
	\vspace*{-1.5em}
	\begin{center}
		\begin{tabular}{c|rr|r}
			\hline\hline
	        Training &\multicolumn{2}{c|}{MSE}&\\\cline{2-3}	        
	        Size& \textbf{P-DNN} & \textbf{SI-DNN} & $\rho$\\
	        \hline\hline
	        100& 2.80 & 1.50 & 10  \\
	        \hline
	        1000& 1.00 & 0.59 & 2   \\
	        \hline
	        5000& 0.54& 0.32 & 1   \\
			\hline
			10000& 0.19& 0.14 & 0.2   \\
			\hline\hline
			\end{tabular}
		\end{center}
		\label{tbl:39large}
		\vspace*{-0.5em}
\end{table}

The performance of the two DNNs was evaluated in terms of the MSE for training sizes $(100,1000,5000,10000)$ sampled from the complete training set of size 10000. The batch-size for all tests was fixed to 100. The cross-validation set was used to determine the initial learning rate (ILR), epochs needed, and the factor $\rho$ in \eqref{eq:si-dnn}. The ILR for training sizes $(100,1000,5000,10000)$ was $(5,5,10,50)\times 10^{-4}$ and the epochs needed were $(2000,500,500,250)$ for both DNNs. The decrease in the training epochs needed is due to the increase in gradient steps per epoch for larger training sizes with fixed batch size. The MSEs obtained by the two DNNs averaged over the 1000 test instances are provided in Table~\ref{tbl:39large} alongside the factor $\rho$ used for different training sizes. As anticipated, the test errors for both DNNs decrease for larger training sizes. However, the SI-DNN consistently outperforms the P-DNN with the improvement being more pronounced at smaller training sizes. Interestingly, a decreasing trend in the suitable choice of $\rho$ was obtained from cross-validation indicating that as the training samples become abundant, sensitivity information seems to be becoming less important. The remaining numerical tests explicitly focus on small training sizes. For simplicity, hereon we fix the ILR to $5\times10^{-4}$ and $\rho=20$.

\subsubsection{Tests on IEEE 39-bus system with small training datasets}

\begin{figure}[t]
    \centering
    \hspace*{-1em}\includegraphics[scale=0.27]{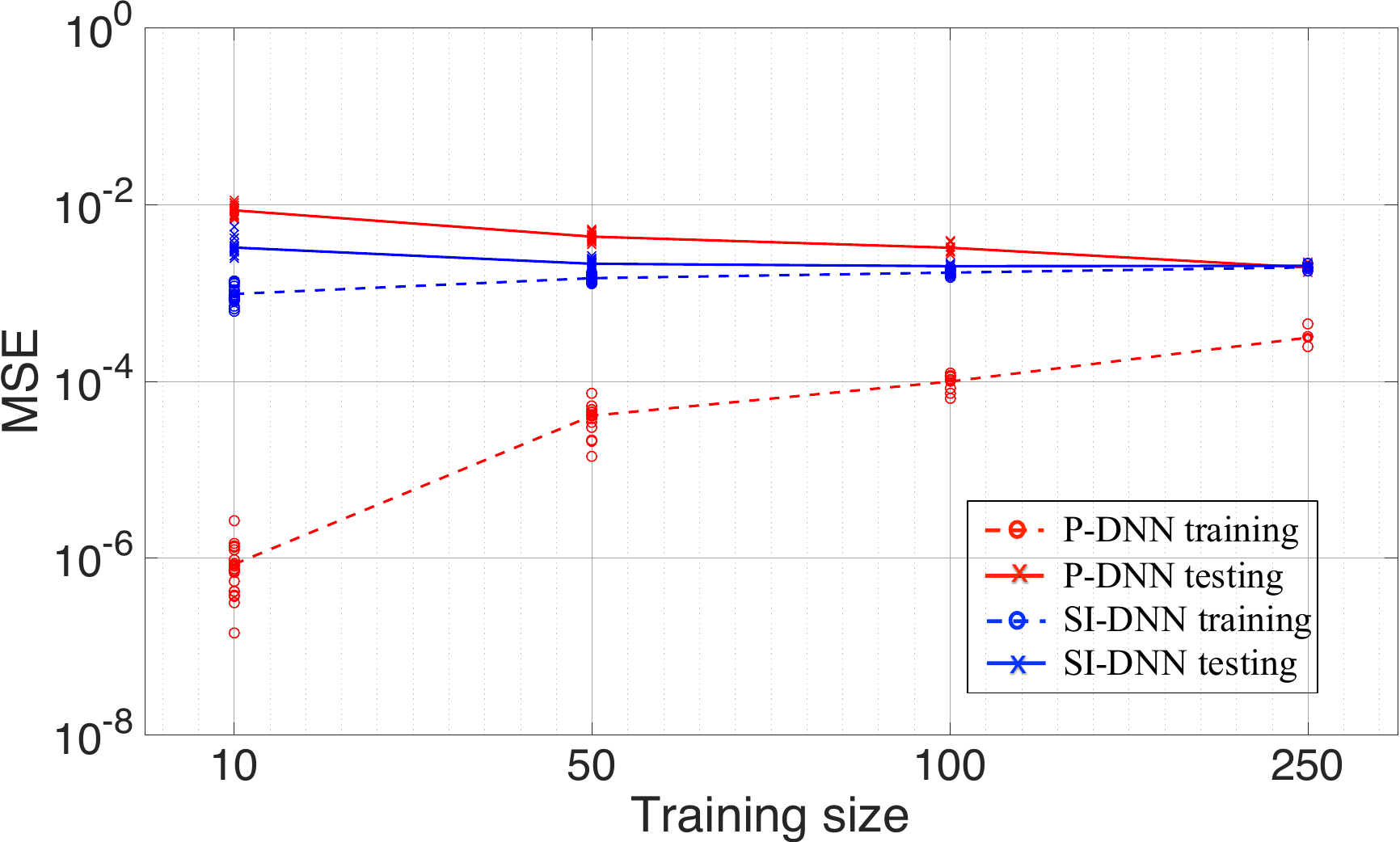}\vspace*{0.4em}
    \includegraphics[scale=0.27]{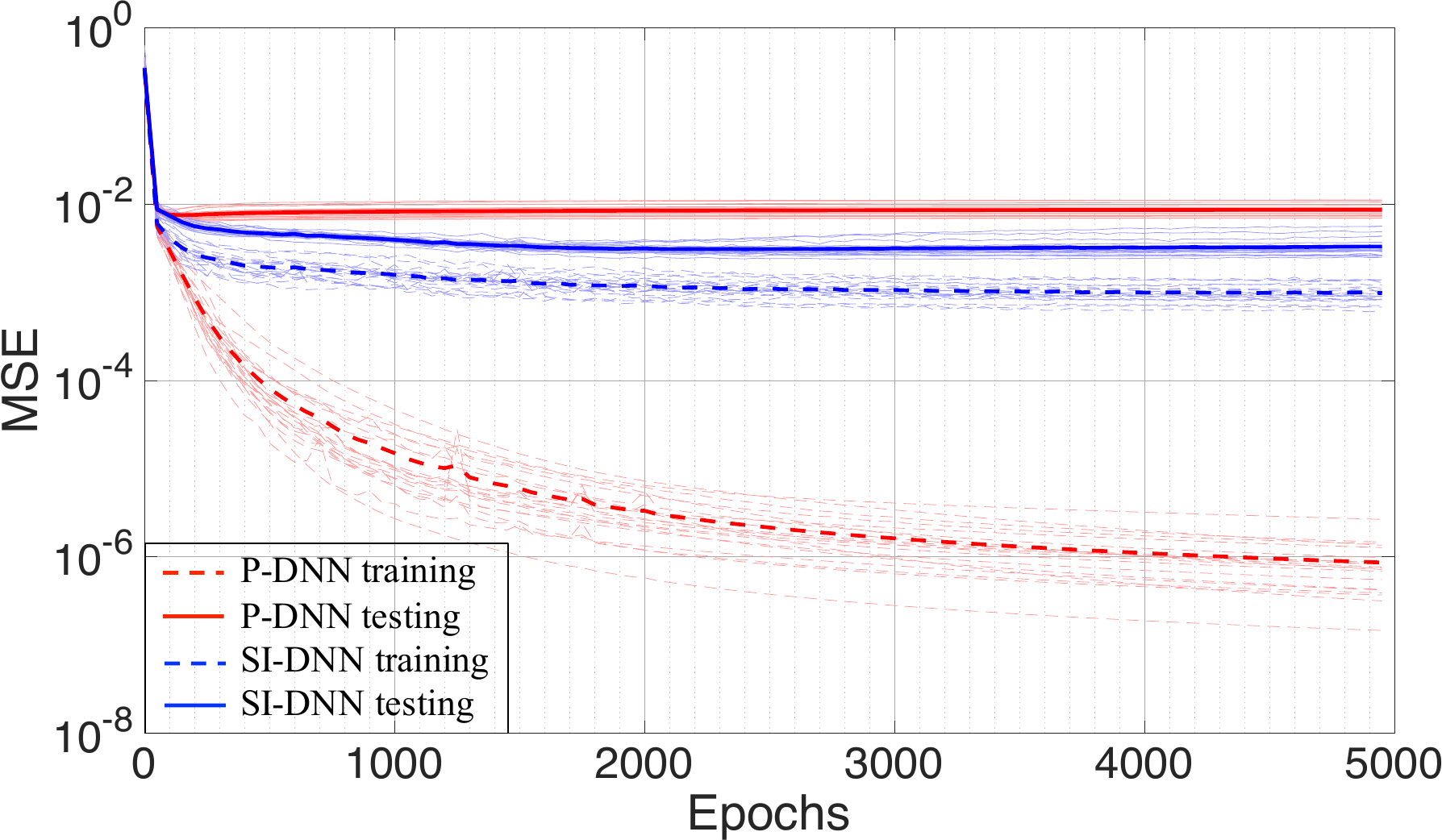}
    \vspace*{-1em}
    \caption{Average training and testing errors for different training sizes (\emph{top}); and errors across epochs for different runs with training size 10~(\emph{bottom}).}
    \label{fig:MP39sum}
\end{figure}

\begin{table}[t]
	\caption{Average Test MSE [$\times~10^{-3}$] and training Time [in sec] for predicting MATPOWER solution on IEEE 39-bus system}
	\vspace*{-1.5em}
	\begin{center}
		\begin{tabular}{c|rr|rr}
			\hline\hline
	        Training &\multicolumn{2}{c|}{\textbf{P-DNN}}&\multicolumn{2}{c}{\textbf{SI-DNN}}\\\cline{2-5}	        
	        Size& MSE & Time & MSE & Time\\
	        \hline\hline
	        10& 8.6 & 738 & 3.3 & 746 \\
	        \hline
	        50& 4.3& 739 & 2.1 & 756  \\
	        \hline
	        100& 3.2& 747 & 2.0 & 776  \\
			\hline
			250& 1.9& 302 & 2.0 & 332  \\
			\hline\hline
			\end{tabular}
		\end{center}
		\label{tbl:MP39}
		\vspace*{-0.5em}
\end{table}

A dataset $\{(\btheta_s,\bJ_{\btheta_s},\cbx_{\btheta_s})\}_{s=1}^{1000}$ was created following the methodology delineated in the previous subsection to evaluate the two DNNs when the training sizes are varied over a data-scarce regime. The evaluation was performed as follows. First, for a \emph{training size} of 10, we created 20 different training sets by sampling 10 OPF instances from the dataset without replacement. For each of these 20 times or \emph{runs}, the OPF instances not sampled for training consisted the testing sets. We then separately trained P-DNN and SI-DNN on these 20 sets. For the training sizes of $(50,100,250)$, we had $(20,10,4)$ runs, respectively. For training sizes $(10, 50, 100)$, the entire training set was used for gradient computation at each step, with the total epochs being 5000. When the training size was 250, the batch-size was fixed to 100, and total epochs to 2000. The  training and testing MSE loss for all training sizes, and runs are shown in Fig.~\ref{fig:MP39sum}~\emph{(top)}. For the tests with training size 10, the evolution of DNN errors are shown in Fig.~\ref{fig:MP39sum}~\emph{(bottom)}. The average test MSE and training times for the two DNNs are shown in Table~\ref{tbl:MP39}. From Fig.~\ref{fig:MP39sum}~\emph{(top)}, we observe as anticipated, that for both DNNs, the gap between training and testing loss decreases for larger training size. Further, the errors for different \emph{runs} are well clustered, indicating a numerically stable DNN implementation. From Table~\ref{tbl:MP39}, it is fascinating to note that the test loss attained by SI-DNN is much lower than P-DNN, especially at smaller training sets. For instance, the P-DNN requires 100 samples to roughly attain the average test MSE which the SI-DNN attains with 10 samples. The lower MSE for P-DNN with training size 250 is a repercussion of not updating $\rho$ for varying training sizes, which was avoided for simplicity. It is worth stressing that the improvement in sample efficiency comes at modest increase in training time.

\subsubsection{Tests on other benchmarks}

\begin{table}[t]
	\caption{Average Test MSE [$\times~10^{-3}$] and Training Time [in sec] for predicting MATPOWER solution}
	\vspace*{-2em}
	\begin{center}
		\begin{tabular}{c|rr|rr|rr|rr}
			\hline\hline
			\multirow{3}{*}{\rotatebox{90}{\parbox{3em}{Train.\\ Size}}}
			&\multicolumn{4}{c|}{\textbf{IEEE 118-bus}}&\multicolumn{4}{c}{\textbf{Illinois 200-bus}}\\ \cline{2-9}
	        &\multicolumn{2}{c|}{{\textbf{P-DNN}}}&\multicolumn{2}{c|}{{\textbf{SI-DNN}}}&\multicolumn{2}{c|}{\textbf{P-DNN}}&\multicolumn{2}{c}{\textbf{SI-DNN}}\\\cline{2-9}	        
	        & MSE & Time& MSE & Time& MSE & Time & MSE & Time\\
	        \hline\hline
	        25& 1.8 &447 &1.1 &483 &0.19 &452 &0.04 &491 \\
	        \hline
	        50& 1.7 &458 &1.1 &527 &0.15 &456 &0.04 &524 \\
	        \hline
			100& 1.6 &463 &0.9 &610 &0.09 &471 &0.06 &608 \\
			\hline\hline
			\end{tabular}
		\end{center}
		\label{tbl:MPall}
		\vspace*{-2em}
\end{table}

The DNN architecture chosen for the other two power systems was similar to the IEEE 39-bus case with the differences being in the number of neurons per layer. Specifically, the DNNs used for the 118- and 200-bus systems had 512 neurons in hidden layers, with the input (output) layers having 128 (107), and 302 (97) neurons, respectively. For each of these systems, we created a dataset with 500 feasible\footnote{To obtain a dataset of 500 instances, the OPF was solved for 550 instances and the first 500 feasible instances were retained. For the 118-bus system, all instances were feasible while for the 200-bus system, four infeasible instances were encountered.} random demands generated by scaling the nominal demands entry-wise by factors drawn uniformly from [0.7,1.3]. The linear cost coefficients from the respective benchmark systems were retained as $c_p$'s, while the reactive power cost coefficients were set to zero. All DNNs were evaluated for five runs, with training sizes of 25, 50, and 100. Table~\ref{tbl:MPall} summarizes the obtained results.

\begin{table}[t]
	\caption{Average time [in sec] to solve AC-OPF, $t_\text{OPF}$; compute sensitivities, $t_\text{sa}$; and obtain DNN predictions followed by running AC power flow, $t_\text{infer}$}
	\vspace*{-1.5em}
	\begin{center}
		\begin{tabular}{c|rrr}
			\hline\hline
	        Test System & $t_\text{OPF}$ & $t_\text{sa}$ & $t_\text{infer}$\\\cline{2-4}
	        \hline
	        39-bus& 0.1229 & 0.0034 & 0.0039 \\
	        \hline
	        118-bus& 0.2577 & 0.0260 & 0.0050  \\
	        \hline
	        200-bus& 0.3032 & 0.0811 & 0.0078   \\
			\hline\hline
			\end{tabular}
		\end{center}
		\label{tbl:times}
		\vspace*{-1.5em}
\end{table}

Having evaluated the improvement in MSE brought by the sensitivity-informed learning approach, we next assessed the additional time-complexity introduced for computing the desired sensitivities. Specifically, while building the datasets for the IEEE 39-bus, the IEEE 118-bus, and the Illinois 200-bus system, we computed: \emph{i)} the average time $t_\text{OPF}$ taken by MATPOWER to solve an OPF instance; \emph{ii)} the average time $t_\text{sa}$ required for computing the Jacobian matrix $\bJ_{\btheta_s}$ using~\eqref{eq:diff}\footnote{For improved numerical performance, matrix $\bS$ was stored as a sparse matrix and \eqref{eq:diff} was solved using MATLAB's command~{\texttt{lsqminnorm}}.}; and \emph{iii)} the average time $t_\text{infer}$ needed to obtain a \emph{complete} OPF minimizer using SI-DNN during the inference phase. To do the latter, we summed up the time taken for evaluating SI-DNN predictions and the time needed to evaluate a corresponding AC power flow solution using MATPOWER. It must be noted that evaluating $t_\text{infer}$ is merely to assess an approximate speed-up offered by the DNNs over conventional OPF solvers. It does not constitute a rigorous comparison since neither optimality nor feasibility is guaranteed for DNN predictions. The aforementioned times are reported in Table~\ref{tbl:times}. It is exciting to observe that while the SI-DNN approach can reduce the training size requirement by up to a factor of 10, evaluating sensitivities for training the SI-DNN requires substantially less time than solving an OPF instance. Finally, the average speed-up factor $t_\text{OPF}/t_\text{infer}$ obtained for the 39-, 118-, and 200-bus systems was approximately 34, 63, and 52, respectively. 

\subsection{Learning Globally Optimal OPF Solutions}\label{subsec:global}


\begin{table}[t]
	\caption{Average Test MSE [$\times~10^{-3}$] for predicting SDP solutions, and constraint violation statistics on the IEEE 39-bus system}
	\begin{center}
		\begin{tabular}{c|rrrr|rrrr}
			\hline\hline
	        Train. &\multicolumn{4}{c|}{\textbf{P-DNN}}&\multicolumn{4}{c}{\textbf{SI-DNN}}\\\cline{2-9}	        
	        Size& MSE & \textbf{(a)} & \textbf{(b)}&\textbf{(c)} &MSE & \textbf{(a)}& \textbf{(b)}& \textbf{(c)}\\
	        \hline\hline
	        10& 6.3 & 2.61 & 0.50 & 9.78& 0.91 & 2.52 & 0.37 & 3.35 \\
	        \hline
	        50& 3.6& 2.45 & 0.55 & 7.38& 0.62 & 2.58 & 0.27 & 2.06  \\
	        \hline
	        100& 2.5& 2.59 & 0.53 & 6.87 &0.67 & 2.52 & 0.27 & 1.96   \\
			\hline
			\multicolumn{9}{c}{\textbf{(a)}~\#violations /instance; \textbf{(b)}~max. violation; \textbf{(c)}~mean violation [$\times10^{-4}$]}\\\cline{1-9}
			\hline\hline
			\end{tabular}
		\end{center}
		\label{tbl:SDP39}
\end{table}

The SI-DNN was evaluated towards predicting the minimizer of an SDP relaxation-based OPF solver for the IEEE 39-bus system. A uniform active power cost was used for all generators while the reactive power cost coefficients were set as $c_n^q=0.1c_n^p$. To build a dataset, a set of $1,000$ random $\btheta's$ was sampled as explained earlier. The corresponding $\cbx_{\theta}$'s were obtained by solving~\eqref{eq:SDP} using the MATLAB-based optimization toolbox YALMIP with SDP solver MOSEK~\cite{YALMIP}. For all SDP instances, the ratio of the second largest eigenvalue of matrix $\bV$ to the largest eigenvalue was found to lie in~$[3\cdot10^{-7},1\cdot10^{-4}]$; numerically indicating an exact relaxation. Thus, the eigenvector corresponding to the largest eigenvalue was deemed as the optimal voltage~$\bv$. If instances with inexact relaxation are encountered, they can be omitted from the dataset. As with learning MATPOWER solutions, the sample efficiency of SI-DNN was found superior to P-DNN in learning globally optimal OPF solutions; see Table~\ref{tbl:SDP39} for the average MSE attained during testing. The presented results with local and global OPF solvers demonstrate that SI-DNN yields a dramatic improvement in generalizability. The feasibility statistics included in Table~\ref{tbl:SDP39} are elaborated upon in Section~\ref{subsec:feasible}.

\subsection{DNN Performance Evaluation under a  Time Budget}\label{subsec:time-budget}

\begin{table}[t]
	\caption{Average test MSE $[\times 10^{-4}]$ for predicting MATPOWER solution and training set generation time-budget $t_\text{gen}$}
	\vspace*{-1.5em}
	\begin{center}
		\begin{tabular}{c|r|r|r|r}
			\hline\hline
	        \textbf{Benchmark} & $t_\text{gen}$
	        &\textbf{P-DNN}&\textbf{SI-DNN}&\textbf{Improvement}\\
	        &\textbf{[sec]} & \textbf{MSE} & \textbf{MSE}&\textbf{\%}\\
	        \hline\hline
	        \multirow{4}{*}{\rotatebox{0}{\parbox{3em}{39-bus}}}
	        &1.2&  0.86 & 0.24 &72.1 \\
	        \cline{2-5}
	        &6.2&  0.49 & 0.21 &57.1 \\
	        \cline{2-5}
	        &12.3& 0.33 & 0.19 &40.6\\
	        \cline{2-5}
	        &30.7&  0.17 & 0.12 &29.4 \\
	        \cline{2-5}
	        \hline\hline
	        \multirow{3}{*}{\rotatebox{0}{\parbox{4em}{118-bus}}}
	        &6.45&  17.50 & 10.80 & 38.3 \\
	        \cline{2-5}
	        &12.89& 17.84 & 12.94 & 27.5  \\
	        \cline{2-5}
	        &25.77& 17.01 & 7.82 & 54.1 \\
			\hline\hline
			\multirow{3}{*}{\rotatebox{0}{\parbox{4em}{200-bus}}}
	        &7.58&  1.83 & 0.32 & 82.5 \\
	        \cline{2-5}
	        &15.16& 1.46 & 0.45 & 69.2  \\
	        \cline{2-5}
	        &30.32& 0.91 & 0.57 & 37.4 \\
			\hline\hline
			\end{tabular}
		\end{center}
		\label{tbl:timebudget}
		\vspace*{-1.5em}
\end{table}

 Tables~\ref{tbl:39large}--\ref{tbl:MPall} and~\ref{tbl:SDP39} attest the improved sample efficiency of sensitivity-informed over conventional training under various dataset sizes, benchmark networks, and OPF solvers. This section exemplifies how these results translate to gains in MSE for a fixed time budget. Specifically, P-DNN and SI-DNN were compared when allotted identical times to complete training. This evaluation was carried out for the settings described in Section~\ref{subsec:localtests}. For all tests, the number of epochs for P-DNN was kept fixed as provided in Section~\ref{subsec:localtests}, while the number of epochs for SI-DNN was reduced, so as to match the training time of the P-DNN. Next, a common time budget $t_\text{gen}$ was fixed for creating the training datasets. Based on Table~\ref{tbl:times}, the training sizes for P-DNN and SI-DNN can be approximately computed as $t_\text{gen}/t_\text{OPF}$ and $t_\text{gen}/(t_\text{OPF}+t_\text{sa})$; implying smaller training sets for the SI-DNN. The test MSEs obtained for the aforementioned setup are provided in Table~\ref{tbl:timebudget}. It was observed that with identical time budgets, an SI-DNN yields test MSEs that are 28-83\% less compared to P-DNN.

\subsection{Assessing Feasibility of DNN Predictions}\label{subsec:feasible}

While emphasis has been on MSE, the importance of satisfying constraints cannot be undermined. To this end, we tested the feasibility of SI-DNN OPF predictions using the following metrics. For each of the DNNs, given a test input and the associated DNN prediction, an AC power flow solution was obtained using MATPOWER. For each instance, the inequalities in~\eqref{eq:P2:m} not directly enforced by the $\tanh$ activation were evaluated. These included voltage limits on load buses, line flow limits, generator reactive power limits, and the slack bus active power limits, totalling to 126, 424, and 647 constraints for the 39-, 118-, and 200-bus system, respectively. For suitable scaling, the violations in flows and generation were normalized by the maximum limit. To be specific, a normalized violation of $10^{-3}$ in generator power injection translates to a violation of $0.1\%$ of the maximum power capacity of that generator. Voltage violations were maintained in pu. 

We first evaluated the constraint violations caused by SI-DNN and P-DNN predictions while learning globally optimal OPF solutions obtained from the SDP-based solver. The assessment was carried out on the test instances that remained after sampling training sets of different sizes from the 1,000 random instances [cf.~Section~\ref{subsec:global}]. For different training sizes, Table~\ref{tbl:SDP39} lists: \emph{a)} the average number of violations exceeding a normalized magnitude of $10^{-6}$ per test instance; \emph{b)} the maximum constraint violation observed; and \emph{c)} the violations averaged over all constraints and test instances. Interestingly, while both DNNs incur similar count of violations, SI-DNN reduces the maximum violation by half and the mean violation to less than one third. We further investigated into the specific constraints being violated by SI-DNN predictions. Interestingly, there were just 5 constraints frequently violated. Three of these were minimum reactive power generation, and the remaining were maximum active power of the slack generator and a line flow limit.

\begin{table}[t]
	\caption{($\mathrm{a}$) Average violations per instance; ($\mathrm{b}$) maximum violation; and ($\mathrm{c}$) mean violation $[\times 10^{-4}]$ for predicting MATPOWER solution}
	\vspace*{-1.5em}
	\begin{center}
		\begin{tabular}{c|c|rrr|rrr}
			\hline\hline
	        \textbf{Benchmark} & \textbf{Train.}
	        &\multicolumn{3}{c|}{\textbf{P-DNN}}&\multicolumn{3}{c}{\textbf{SI-DNN}}\\\cline{3-8}	        
	        &\textbf{Size}& \textbf{(a)} & \textbf{(b)}&\textbf{(c)} & \textbf{(a)}& \textbf{(b)}& \textbf{(c)}\\
	        \hline\hline
	        \multirow{4}{*}{\rotatebox{0}{\parbox{3em}{39-bus}}}
	        &10&  2.29 & 1.01 & 15&  1.91 & 0.33 & 7.66 \\
	        \cline{2-8}
	        &50& 2.45 & 0.95 & 9.55& 2.14 & 0.24 & 6.30  \\
	        \cline{2-8}
	        &100& 2.31 & 0.71 & 7.65 & 1.96 & 0.26 & 6.26   \\
	        \cline{2-8}
	        &250& 2.23 & 0.51 & 6.99 & 1.87 & 0.28 & 6.89   \\
	        \hline\hline
	        \multirow{3}{*}{\rotatebox{0}{\parbox{4em}{200-bus}}}
	        &25&  10.43 & 1.42 & 8.10&  4.99 & 1.37 & 4.27 \\
	        \cline{2-8}
	        &50& 9.08 & 1.19 & 6.17& 3.69 & 1.14 & 4.64 \\
	        \cline{2-8}
	        &100& 10.89 & 1.15 & 7.17 & 2.74 & 1.24 & 4.71 \\
			\hline\hline
			\end{tabular}
		\end{center}
		\label{tbl:MPviolations}
		\vspace*{-1.5em}
\end{table}

\begin{table}[t]
	\caption{($\mathrm{a}$) Average violations per instance; ($\mathrm{b}$) maximum violation; and ($\mathrm{c}$) mean violation $[\times 10^{-3}]$ for predicting MATPOWER solution on IEEE 118-bus system}
	\vspace*{-1.5em}
	\begin{center}
		\begin{tabular}{c|c|rrr|rrr}
			\hline\hline
	        \textbf{Constraint} & \textbf{Train.}
	        &\multicolumn{3}{c|}{\textbf{P-DNN}}&\multicolumn{3}{c}{\textbf{SI-DNN}}\\\cline{3-8}	        
	        \textbf{Set}&\textbf{Size}& \textbf{(a)} & \textbf{(b)}&\textbf{(c)} & \textbf{(a)}& \textbf{(b)}& \textbf{(c)}\\
	        \hline\hline
	        \multirow{3}{*}{\rotatebox{0}{\parbox{2em}{Full}}}
	        &25&  7.97 & 1.08 & 1.90&  9.28 & 1.44 & 3.90 \\
	        \cline{2-8}
	        &50& 8.24 & 1.02 & 1.60& 10.20 & 0.97 & 3.30  \\
	        \cline{2-8}
	        &100& 8.45 & 1.18 & 1.80 & 10.09 & 0.93 & 2.20   \\
	        \hline\hline
	        \multirow{3}{*}{\rotatebox{0}{\parbox{4em}{Reduced}}}
	        &25&  2.12 & 0.73 & 0.32&  1.67 & 0.66 & 0.14 \\
	        \cline{2-8}
	        &50& 1.92 & 0.63 & 0.21& 1.65 & 0.68 & 0.14 \\
	        \cline{2-8}
	        &100& 2.01 & 1.17 & 0.34 & 1.77 & 0.78 & 0.12 \\
			\hline\hline
			\end{tabular}
		\end{center}
		\label{tbl:MPviolations118}
		\vspace*{-1.5em}
\end{table}

We repeated the previous feasibility analysis for the DNNs aimed at learning locally optimal OPF solutions from MATPOWER. The constraint violation statistics obtained for the IEEE 39-bus and the Illinois 200-bus system, provided in Table~\ref{tbl:MPviolations}, consistently demonstrate the improvements yielded by the SI-DNN approach for different training sizes. Table~\ref{tbl:MPviolations118} reports the same statistics for the IEEE 118-bus system while predicting the MATPOWER solution. The numerical observations for the IEEE 118-bus system do not align with the results for other benchmark networks. The statistics provided in the top part of Table~\ref{tbl:MPviolations118} exhibit much higher constraint violations for both P-DNN and SI-DNN; note that the mean violations are of the order $10^{-3}$ as opposed to $10^{-4}$ for other networks. Moreover, SI-DNN performs worse than P-DNN on several metrics. Spurred by the exceptionally high constraint violations, we investigated the individual constraints being violated. It was found that several generator reactive power limits were being consistently violated by both P-DNN and SI-DNN. It turns out that these limits were binding for \emph{all} the random OPF scenarios in the training and testing datasets. For benchmarks that exhibit such patterns with certain dispatch quantities being fixed across scenarios, it may be prudent to set them at the respective values and solve a reduced OPF. To emphasize on the violation statistics for the non-trivial constraints, we computed the feasibility metrics on a reduced set of constraints not including those reactive power limits that were consistently binding. The obtained results shown at the bottom of Table~\ref{tbl:MPviolations118} corroborate the superior performance of SI-DNN over P-DNN.

\section{Conclusions}\label{sec:conclusions}
This work has built on the fresh idea of sensitivity-informed training for learning the solutions of arbitrary AC-OPF formulations. It comprehensively delineated the steps for computing the involved sensitivities using the optimal primal/dual solutions, which are readily available by AC-OPF solvers. Such sensitivities of the primal AC-OPF solutions have been shown to exist under mild assumptions, while their computation is as simple as solving a system of linear equations with multiple right-hand sides. The approach is quite general since the OPF solutions comprising the training dataset can be obtained by off-the-shelf nonlinear OPF solvers or modern conic relaxation-based schemes. It is also worth stressing that sensitivity-informed training can readily complement other existing learn-to-OPF methodologies. Extensive numerical tests on three benchmark power systems have demonstrated that with a modest increase in training time, SI-DNNs attain the same prediction performance as conventionally trained DNNs by using roughly only 1/10 to 1/4 of the training data. Such improvement on sample efficiency reduces the time needed for generating training datasets, and is thus, relevant to delay-critical power systems applications. Furthermore, SI-DNN predictions turn out to feature better constraint satisfaction capabilities too. Sensitivity-informed learning forms the solid foundations for several exciting and practically relevant research directions, such as warm-starting key optimal primal/dual variables to accelerate decentralized OPF solvers and predicting active constraints.

\appendix
\subsection{Python Implementation for SI-DNN}\label{sec:AppA}
A typical implementation example for computing the gradient of the MSE loss in P-DNN with respect to the DNN weights (which are the \emph{trainable variables}) involves
\begin{Verbatim}[fontsize=\small]
with tensorflow.GradientTape() as tape:
    pred_x = model(theta)
    loss = keras.losses.MSE(xlabel,pred_x)
model_gradients=tape.gradient(loss,
                model.trainable_variables)
\end{Verbatim}
where {\verb|model|} represents the DNN and {\verb|GradientTape|} computes the desired gradient. In transitioning to SI-DNN, we first need to compute the gradient of the DNN output {\verb|pred_x|} with respect to its input {\verb|theta|} to define the loss. We then compute the gradients of the two loss terms with respect to the DNN weights. This can be implemented using nested {\verb|GradientTape|} as
\begin{Verbatim}[fontsize=\small]
with tensorflow.GradientTape() as tape:
  with tensorflow.GradientTape() as tape2:
    tape2.watch(theta)
    pred_x=model(theta)
    Ploss=keras.losses.MSE(xlabel,pred_x)
  J_model=tape2.batch_jacobian(pred_x,
                                    theta)
  J_flat=tf.keras.backend.reshape(fgrad,
                               shape=(1,))
  SI_loss=keras.losses.MSE(J_flat,J_label)
  total_loss=P_loss+rho*SI_loss
model_gradients=tape.gradient(loss,
                model.trainable_variables)
\end{Verbatim}
where the inner {\verb|tape|} computes the sensitivity of DNN to compute the overall SI-DNN loss, while the outer {\verb|tape|} computes the gradients for weight updates.

\subsection{Sensitivity computation with MATPOWER}\label{sec:AppB}
While solving the AC-OPF instances with MATPOWER, we used the Cartesian coordinate system, and flow limits were imposed on squared currents. For computing the desired sensitivities, we first need to build the linear system~\eqref{eq:diff}, which requires the optimal dual variables, constraint function values, and the derivatives of the constraint functions with respect to the optimization variables. Although MATPOWER can deal with the AC-OPF posed with voltages in Cartesian coordinates, it slightly differs from the QCQP in~\eqref{eq:P1} as follows:
\renewcommand{\labelenumi}{\emph{a\arabic{enumi})}}
\begin{enumerate}
    \item MATPOWER enforces power flow equations as in \eqref{eq:P1:pg}--\eqref{eq:P1:ql}. Different from \eqref{eq:P1:pglim}--\eqref{eq:P1:qglim} however, MATPOWER poses generator (re)active power limits as $\underline{p}_n^g\leq p_n^g\leq \bar{p}_n^g$. Thus, the related $\nabla_\bv g_m$ becomes zero and $\nabla_{\bx_g} g_m$ becomes a signed canonical vector corresponding to bus $n$.
    \item MATPOWER constraints voltages, rather than squared voltages as in~\eqref{eq:P1:vlim}. While the two versions are equivalent, the related derivatives $\nabla_\bv g$ apparently differ. The derivatives of non-squared magnitudes can be found using the chain rule as $\nabla_\bv v_n=\frac{\partial v_n}{\partial v_n^2}\nabla_\bv v_n^2= \frac{1}{v_n}\nabla_\bv v_n^2$ and $\nabla_\bv v_n^2=2\bM_{v_n}\bv$ from \eqref{eq:volts}.
    \item Different from~\eqref{eq:P1:ref}, MATPOWER sets the voltage angle reference to zero by enforcing $\arctan(v_{N+1}^i/v_{N+1}^r)=0$. Fortunately, simply setting the imaginary part $v_{N+1}^i$ to zero is equivalent, and the gradients of these two formulations agree. Thus, despite the difference in formulation, we use~\eqref{eq:P1:ref} wherever needed in building~\eqref{eq:diff}.
    \item Finally, MATPOWER poses flow limits on both the sending and receiving ends of each line, thus doubling the number of constraints in~\eqref{eq:P1:flim}. The matrices $\bM_{\bi_{mn}}$ can be built using the \emph{from bus} and \emph{to bus} admittances obtained via the MATPOWER command~{\verb|makeYbus()|}.
\end{enumerate}

\balance
\bibliography{myabrv,power,inverters}

\begin{thebibliography}{10}
\providecommand{\url}[1]{#1}
\csname url@samestyle\endcsname
\providecommand{\newblock}{\relax}
\providecommand{\bibinfo}[2]{#2}
\providecommand{\BIBentrySTDinterwordspacing}{\spaceskip=0pt\relax}
\providecommand{\BIBentryALTinterwordstretchfactor}{4}
\providecommand{\BIBentryALTinterwordspacing}{\spaceskip=\fontdimen2\font plus
\BIBentryALTinterwordstretchfactor\fontdimen3\font minus
  \fontdimen4\font\relax}
\providecommand{\BIBforeignlanguage}[2]{{%
\expandafter\ifx\csname l@#1\endcsname\relax
\typeout{** WARNING: IEEEtran.bst: No hyphenation pattern has been}%
\typeout{** loaded for the language `#1'. Using the pattern for}%
\typeout{** the default language instead.}%
\else
\language=\csname l@#1\endcsname
\fi
#2}}
\providecommand{\BIBdecl}{\relax}
\BIBdecl

\bibitem{MATPOWER}
R.~D. Zimmerman, C.~E. Murillo-Sanchez, and R.~J. Thomas, ``{MATPOWER}:
  steady-state operations, planning and analysis tools for power systems
  research and education,'' \emph{{IEEE} Trans. Power Syst.}, vol.~26, no.~1,
  pp. 12--19, Feb. 2011.

\bibitem{Bai08}
X.~Bai, H.~Wei, K.~Fujisawa, and Y.~Yang, ``Semidefinite programming for
  optimal power flow problems,'' \emph{{Intl. Journal of Electric Power \&
  Energy Systems}}, vol.~30, no.~6, pp. 383--392, 2008.

\bibitem{Bose}
S.~Bose, D.~Gayme, K.~Chandy, and S.~Low, ``Quadratically constrained quadratic
  programs on acyclic graphs with application to power flow,'' \emph{{IEEE}
  Trans. Control of Network Systems}, vol.~2, no.~3, pp. 278--287, Sep. 2015.

\bibitem{MSL15}
R.~{Madani}, S.~{Sojoudi}, and J.~{Lavaei}, ``Convex relaxation for optimal
  power flow problem: Mesh networks,'' \emph{{IEEE} Trans. Power Syst.},
  vol.~30, no.~1, pp. 199--211, Jan. 2015.

\bibitem{Low14}
S.~Low, ``Convex relaxation of optimal power flow -- {Part II: Exactness},''
  \emph{{IEEE} Trans. Control of Network Systems}, vol.~1, no.~2, pp. 177--189,
  Jun. 2014.

\bibitem{Xavier20SCUC}
A.~S. Xavier, F.~Qiu, and S.~Ahmed, ``Learning to solve large-scale
  security-constrained unit commitment problems,'' \emph{INFORMS Journal on
  Computing}, pp. 1--18, Oct. 2020, (early access).

\bibitem{Fioretto1}
F.~Fioretto, T.~W. Mak, and P.~V. Hentenryck, ``Predicting {AC} optimal power
  flows: Combining deep learning and {Lagrangian} dual methods,'' in
  \emph{{AAAI Conf. on Artificial Intelligence}}, New York, NY, Feb. 2020.

\bibitem{DeepOPFPan19}
X.~Pan, T.~Zhao, and M.~Chen, ``{DeepOPF}: Deep neural network for {DC} optimal
  power flow,'' in \emph{Proc. {IEEE} Intl. Conf. on Smart Grid Commun.},
  Beijing, China, Oct. 2019, pp. 1--6.

\bibitem{Pan20DeepOPFplus}
T.~Zhao, X.~Pan, M.~Chen, A.~Venzke, and S.~H. Low, ``Deepopf+: A deep neural
  network approach for {DC} optimal power flow for ensuring feasibility,'' in
  \emph{Proc. {IEEE} Intl. Conf. on Smart Grid Commun.}, Tempe, AZ, Nov. 2020,
  pp. 1--6.

\bibitem{pan2020feasopt}
\BIBentryALTinterwordspacing
X.~Pan, M.~Chen, T.~Zhao, and S.~H. Low, ``Deepopf: A feasibility-optimized
  deep neural network approach for {AC} optimal power flow problems,'' 2020,
  (preprint). [Online]. Available: \url{https://arxiv.org/abs/2007.01002}
\BIBentrySTDinterwordspacing

\bibitem{GuhaACOPF}
\BIBentryALTinterwordspacing
N.~Guha, Z.~Wang, M.~Wytock, and A.~Majumdar, ``Machine learning for {AC}
  optimal power flow,'' 2019, climate Change Workshop at ICML 2019. [Online].
  Available: \url{https://arxiv.org/abs/1910.08842}
\BIBentrySTDinterwordspacing

\bibitem{ZamzamBaker19}
A.~Zamzam and K.~Baker, ``Learning optimal solutions for extremely fast {AC}
  optimal power flow,'' in \emph{Proc. {IEEE} Intl. Conf. on Smart Grid
  Commun.}, Tempe, AZ, Nov. 2020, pp. 1--6.

\bibitem{YueZhao_DL4P}
Y.~Zhao and B.~Zhang, ``Deep learning in power systems,'' in \emph{Advanced
  Data Analytics for Power Systems}, A.~Tajer, S.~M. Perlaza, and H.~V. Poor,
  Eds.\hskip 1em plus 0.5em minus 0.4em\relax Cambridge, UK: Cambridge
  University Press, May 2021.

\bibitem{RibeiroGNNOPF}
D.~Owerko, F.~Gama, and A.~Ribeiro, ``Optimal power flow using graph neural
  networks,'' in \emph{Proc. {IEEE} Intl. Conf. on Acoustics, Speech, and
  Signal Process.}, Barcelona, Spain, May 2020, pp. 5930--5934.

\bibitem{GKJ2020}
S.~Gupta, V.~Kekatos, and M.~Jin, ``Communication-limited inverter control
  using deep neural networks,'' in \emph{Proc. {IEEE} Intl. Conf. on Smart Grid
  Commun.}, Tempe, AZ, Nov. 2020, pp. 1--6.

\bibitem{lange2020learning}
\BIBentryALTinterwordspacing
H.~Lange, B.~Chen, M.~Berges, and S.~Kar, ``Learning to solve {AC} optimal
  power flow by differentiating through holomorphic embeddings,'' 2020,
  (submitted). [Online]. Available: \url{https://arxiv.org/abs/2012.096224}
\BIBentrySTDinterwordspacing

\bibitem{GMDK21}
\BIBentryALTinterwordspacing
S.~Gupta, S.~Misra, D.~Deka, and V.~Kekatos, ``{DNN}-based policies for
  stochastic {AC-OPF},'' in \emph{Proc. Power Syst. Comput. Conf.}, Porto,
  Portugal, Jun. 2021, (to appear also in the Elsevier Electric Power Systems
  Research). [Online]. Available:
  \url{https://www.faculty.ece.vt.edu/kekatos/papers/PSCC2022a.pdf}
\BIBentrySTDinterwordspacing

\bibitem{OPFandLearnTSG21}
\BIBentryALTinterwordspacing
S.~Gupta, V.~Kekatos, and M.~Jin, ``Controlling smart inverters using proxies{:
  A} chance-constrained {DNN}-based approach,'' \emph{{IEEE} Trans. Smart
  Grid}, May 2021, (submitted). [Online]. Available:
  \url{https://arxiv.org/abs/2105.00429}
\BIBentrySTDinterwordspacing

\bibitem{chen2020metalearning}
\BIBentryALTinterwordspacing
Y.~Chen, S.~Lakshminarayana, C.~Maple, and H.~V. Poor, ``A meta-learning
  approach to the optimal power flow problem under topology reconfigurations,''
  2020, (preprint). [Online]. Available: \url{https://arxiv.org/abs/2012.11524}
\BIBentrySTDinterwordspacing

\bibitem{ChenZhang20}
\BIBentryALTinterwordspacing
Y.~Chen and B.~Zhang, ``Learning to solve network flow problems via neural
  decoding,'' 2020, preprint. [Online]. Available:
  \url{https://arxiv.org/abs/2002.04091}
\BIBentrySTDinterwordspacing

\bibitem{DekaMisraPowerTech19}
D.~{Deka} and S.~{Misra}, ``Learning for {DC-OPF}: {C}lassifying active sets
  using neural nets,'' in \emph{{IEEE PowerTech}}, Milan, Italy, Jun. 2019, pp.
  1--6.

\bibitem{Nandwani19}
M.~Yatin~Nandwani, Abhishek~Pathak and P.~Singla, ``A primal dual formulation
  for deep learning with constraints,'' in \emph{Proc. of Adv. Neural Inf.
  Process. Syst.}, Vancouver, Canada, Dec. 2019, pp. 12\,157--12\,168.

\bibitem{Zhang20icnn}
\BIBentryALTinterwordspacing
L.~Zhang, Y.~Chen, and B.~Zhang, ``A convex neural network solver for {DCOPF}
  with generalization guarantees,'' 2020, (submitted). [Online]. Available:
  \url{https://arxiv.org/abs/2009.09109}
\BIBentrySTDinterwordspacing

\bibitem{ZhangWangGiannakis19}
L.~{Zhang}, G.~{Wang}, and G.~B. {Giannakis}, ``Real-time power system state
  estimation and forecasting via deep unrolled neural networks,'' \emph{{IEEE}
  Trans. Signal Processing}, vol.~67, no.~15, pp. 4069--4077, Aug. 2019.

\bibitem{yang2020robust}
\BIBentryALTinterwordspacing
Q.~Yang, A.~Sadeghi, G.~Wang, G.~B. Giannakis, and J.~Sun, ``Robust {PSSE}
  using graph neural networks for data-driven and topology-aware priors,''
  2020, (submitted). [Online]. Available:
  \url{https://arxiv.org/abs/2003.01667}
\BIBentrySTDinterwordspacing

\bibitem{Raissi19}
M.~Raissi, P.~Perdikaris, and G.~E. Karniadakis, ``Physics-informed neural
  networks: {A} deep learning framework for solving forward and inverse
  problems involving nonlinear partial differential equations,'' \emph{Journal
  of Computational Physics}, vol. 378, pp. 686--707, 2019.

\bibitem{misyris2020PINN}
G.~S. Misyris, A.~Venzke, and S.~Chatzivasileiadis, ``Physics-informed neural
  networks for power systems,'' in \emph{Proc. {IEEE} PES General Meeting},
  Montreal, Canada, Aug. 2020, pp. 1--5.

\bibitem{SGKCB2020}
M.~K. Singh, S.~Gupta, V.~Kekatos, G.~Cavraro, and A.~Bernstein, ``Learning to
  optimize power distribution grids using sensitivity-informed deep neural
  networks,'' in \emph{Proc. {IEEE} Intl. Conf. on Smart Grid Commun.}, Tempe,
  AZ, Nov. 2020, pp. 1--6.

\bibitem{shapiro2013perturbation}
J.~F. Bonnans and A.~Shapiro, \emph{Perturbation Analysis of Optimization
  Problems}.\hskip 1em plus 0.5em minus 0.4em\relax New York, NY: Springer
  Science \& Business Media, 2000.

\bibitem{fiacco1976sensitivity}
A.~V. Fiacco, ``Sensitivity analysis for nonlinear programming using penalty
  methods,'' \emph{Mathematical Programming}, vol.~10, no.~1, pp. 287--311,
  Dec. 1976.

\bibitem{Conejo06}
A.~J. Conejo, E.~Castillo, R.~Minguez, and R.~Garcia-Bertrand,
  \emph{Decomposition Techniques in Mathematical Programming}.\hskip 1em plus
  0.5em minus 0.4em\relax Springer, 2006.

\bibitem{Sidiropoulos18}
H.~Sun, X.~Chen, Q.~Shi, M.~Hong, X.~Fu, and N.~D. Sidiropoulos, ``Learning to
  optimize: {T}raining deep neural networks for interference management,''
  \emph{{IEEE} Trans. Signal Processing}, vol.~66, no.~20, pp. 5438--5453, Oct.
  2018.

\bibitem{BBM03}
F.~Borrelli, A.~Bemporad, and M.~Morari, ``Geometric algorithm for
  multiparametric linear programming,'' \emph{Journal of Optimization Theory
  and Applications}, vol. 118, no.~3, pp. 515--540, Sep. 2003.

\bibitem{TJKT20}
S.~Taheri, M.~Jalali, V.~Kekatos, and L.~Tong, ``Fast probabilistic hosting
  capacity analysis for active distribution systems,'' \emph{{IEEE} Trans.
  Smart Grid}, 2020, (early access).

\bibitem{autodiff}
A.~G. Baydin, B.~A. Pearlmutter, A.~A. Radul, and J.~M. Siskind, ``Automatic
  differentiation in machine learning: {A} survey,'' \emph{{J. Mach. Learn.
  Res.}}, vol.~18, no.~1, pp. 5595--5637, Jan. 2017.

\bibitem{redux}
V.~Kekatos, G.~Wang, H.~Zhu, and G.~B. Giannakis, ``{PSSE} redux: {C}onvex
  relaxation, decentralized, robust, and dynamic approaches,'' in
  \emph{Advances in Power System State Estimation}, M.~El-Hawary, Ed.\hskip 1em
  plus 0.5em minus 0.4em\relax Wiley, 2021.

\bibitem{Almeida94parametric}
K.~Almeida, F.~Galiana, and S.~Soares, ``A general parametric optimal power
  flow,'' \emph{{IEEE} Trans. Power Syst.}, vol.~9, no.~1, pp. 540--547, Feb.
  1994.

\bibitem{Ajjarapu95OCPF}
V.~Ajjarapu and N.~Jain, ``Optimal continuation power flow,'' \emph{Electric
  Power Systems Research}, vol.~35, no.~1, pp. 17--24, Oct. 1995.

\bibitem{Almeida2000varload}
K.~Almeida and R.~Salgado, ``Optimal power flow solutions under variable load
  conditions,'' \emph{{IEEE} Trans. Power Syst.}, vol.~15, no.~4, pp.
  1204--1211, Nov. 2000.

\bibitem{Castillo06}
E.~Castillo, A.~J. Conejo, C.~Castillo, R.~Minguez, and D.~Ortigosa,
  ``Perturbation approach to sensitivity analysis in mathematical
  programming,'' \emph{{Journal of Optimization Theory and Applications}}, vol.
  128, no.~1, pp. 49--74, Jan. 2006.

\bibitem{agrawal2020differentiating}
\BIBentryALTinterwordspacing
A.~Agrawal, S.~Barratt, S.~Boyd, E.~Busseti, and W.~M. Moursi,
  ``Differentiating through a cone program,'' 2020, (submitted). [Online].
  Available: \url{https://arxiv.org/abs/1904.09043}
\BIBentrySTDinterwordspacing

\bibitem{AmosKotler17}
B.~Amos and J.~Z. Kolter, ``{OptNet}: Differentiable optimization as a layer in
  neural networks,'' in \emph{Intl. Conf. on Machine Learning}, Sydney, NSW,
  Australia, 2017, p. 136–145.

\bibitem{Be99}
D.~P. Bertsekas, \emph{Nonlinear Programming}, 2nd~ed.\hskip 1em plus 0.5em
  minus 0.4em\relax Belmont, MA: Athena Scientific, 1999.

\bibitem{Almeida16ill}
K.~C. Almeida and A.~Kocholik, ``Solving ill-posed optimal power flow problems
  via {Fritz-John} optimality conditions,'' \emph{{IEEE} Trans. Power Syst.},
  vol.~31, no.~6, pp. 4913--4922, Nov. 2016.

\bibitem{Dorfler18LICQ}
A.~Hauswirth, S.~Bolognani, G.~Hug, and F.~Dorfler, ``Generic existence of
  unique lagrange multipliers in {AC} optimal power flow,'' \emph{{IEEE} Contr.
  Syst. Lett.}, vol.~2, no.~4, pp. 791--796, Oct. 2018.

\bibitem{Molzhan14}
D.~K. Molzahn, B.~C. Lesieutre, and C.~L. DeMarco, ``A sufficient condition for
  global optimality of solutions to the optimal power flow problem,''
  \emph{{IEEE} Trans. Power Syst.}, vol.~29, no.~2, pp. 978--979, Mar. 2014.

\bibitem{YALMIP}
J.~Lofberg, ``Yalmip : A toolbox for modeling and optimization in matlab,'' in
  \emph{Proc. of the CACSD Conf.}, Taipei, Taiwan, 2004.

\end{thebibliography}
\bibliographystyle{IEEEtran}

\end{document}